\documentclass[leqno]{article}
\usepackage{hyperref}
\usepackage{amssymb,amsmath}
\usepackage{enumerate}

\usepackage[dvipdfmx]{graphicx}
\usepackage{graphicx}
\usepackage{tikz}

\newtheorem{thm}{Theorem}[section]
\newtheorem{cor}[thm]{Corollary}
\newtheorem{lem}[thm]{Lemma}

\newtheorem{defin}[thm]{Definition}
\newtheorem{rem}[thm]{Remark}

\numberwithin{equation}{section}  

\begin{document}

\date{}

\title{On the Cauchy problem for differential operators with double characteristics, a transition from effective to non-effective characteristics}

\author{Tatsuo  Nishitani}

\maketitle

\begin{abstract}
We discuss the well-posedness of the Cauchy problem for hyperbolic operators with double characteristics which changes from non-effectively hyperbolic to effectively hyperbolic, on the double characterisitic manifold,  across a submanifold of codimension $1$. We assume that there is no bicharacteristic tangent to the double characteristic manifold and the spatial dimension is $2$. Then we prove the well-posedness of the Cauchy problem in all Gevrey classes assuming, on the double characteristic manifold, that the ratio of the imaginary part of the subprincipal symbol to the real eigenvalue of the Hamilton map is bounded and that the sum of the real part of the subprincipal symbol and the modulus of the imaginary eigenvalue of the Hamilton map  is strictly positive.
\end{abstract}

\def\dif{\partial}
\def\al{\alpha}
\def\be{\beta}
\def\ga{\gamma}
\def\Ga{\Gamma}
\def\Si{\Sigma}
\def\te{\theta}
\def\ep{\epsilon}
\def\de{\delta}
\def\varep{\varepsilon}
\def\hte{{\hat\theta}}
\def\xii{\langle{\xi'}\rangle}
\def\mxi{\langle{\mu\xi'}\rangle}
\def\xim{\langle{\xi'}\rangle_{\mu}}
\def\mD{\langle{\mu D'}\rangle}
\def\Dm{\langle{D'}\rangle_{\mu}}
\def\bg{{\bar g}}
\def\lr#1{\langle{#1}\rangle}
\def\R{{\mathbb R}}
\def\C{{\mathbb C}}
\def\N{{\mathbb N}}
\def\brho{{\bar\rho}}
\def\la{\lambda}
\def\lam{\lambda_{\mu}}
\def\La{\Lambda}
\def\tLa{{\tilde\Lambda}}
\def\tQ{{\tilde Q}}
\def\qed{\hfill{$\Box$}\vspace{2mm}}
\def\Ker{{\rm Ker}}
\def\hphi{{\hat \phi}}
\def\hg{{\hat g}}
\def\om{\omega}

\section{Introduction}

This paper is a continuation of our  previous papers \cite{Ni3, Ni4}.  Let 
\[
P(x,D)=-D_0^2+A_1(x,D')D_0+A_2(x,D')
\]
be a differential operator of order $2$ in $D_0$ with coefficients $A_j(x,D')$, classical pseudodifferential operator of order $j$ on $\R^n$ depending smoothly on $x_0$  where  $x=(x_0,x')=(x_0,x_1,\ldots,x_n)$. We assume that the principal symbol $p(x,\xi)$  of $P(x,D)$ vanishes exactly of order $2$ on a $C^{\infty}$ manifold $\Sigma$ and
\begin{equation}
\label{eq:rank}
{\rm rank}\,\big(\sum_{j=0}^n d\xi_j\wedge dx_j\big|_{\Sigma}\big)={\rm constant}.
\end{equation}
As in \cite{Ni3, Ni4} we assume that  ${\rm codim}\,\Sigma=3$ and 
\begin{equation}
\label{eq:assump}
\left\{\begin{array}{ll}
\mbox{the spectral structure of $F_p$ changes simply}\\  
\mbox{across a submanifold $S$ of codimension $1$ of $\Sigma$}.
\end{array}\right.
\end{equation}
By conjugation with a Fourier integral operator one can assume $A_1=0$ then, near any point $\rho\in \Sigma$, one can write
\[
p(x,\xi)=-\xi_0^2+\phi_1(x,\xi')^2+\phi_2(x,\xi')^2
\]
where $d\phi_1$ and $d\phi_2$ are linearly independent at $\rho$ and $\Sigma=\{\xi_0=0,\phi_1=0,\phi_2=0\}$. 
Under the assumptions \eqref{eq:rank} and \eqref{eq:assump} without restrictions we can assume (see \cite{Ni3})
\[
\{\xi_0,\phi_2\}>0,\quad \{\xi_0,\phi_1\}=O(|\phi|)
\]
near $\rho$. Here and in what follows $f=O(|\phi|)$, $\phi=(\phi_1,\phi_2)$ means that $f$ is a linear combination of $\phi_1$ and $\phi_2$ near the reference point. We first recall
\begin{lem}
\label{lem:FfF}
{\rm (\cite[Lemma 1.2]{Ni4})}
If the spectral structure of $F_p$ changes across $S$ then we have $\{\xi_0,\phi_2\}^2-\{\phi_1,\phi_2\}^2=0$ on $S$ and 
one of the following cases occurs;
\begin{description}
\item{\rm(i)}
$\{\xi_0,\phi_2\}^2-\{\phi_1,\phi_2\}^2<0$ in $\Sigma\setminus S$ so that $p$ is non-effectively hyperbolic in $\Sigma$  with ${\rm Ker}F_p^2\cap {\rm Im}F_p^2=\{0\}$ in $\Sigma\setminus S$ and ${\rm Ker}F_p^2\cap {\rm Im}F_p^2\neq \{0\}$ on $S$,
\item{\rm(ii)} $\{\xi_0,\phi_2\}^2-\{\phi_1,\phi_2\}^2>0$ in $\Sigma\setminus S$ so that $p$ is effectively hyperbolic in $\Sigma\setminus S$  and non-effectively hyperbolic on $S$ with ${\rm Ker}F_p^2\cap {\rm Im}F_p^2\neq \{0\}$,
\item{\rm(iii)} $\{\xi_0,\phi_2\}^2-\{\phi_1,\phi_2\}^2$ changes the sign across $S$, that is $p$ is effectively hyperbolic in the one side of $\Sigma\setminus S$, non-effectively hyperbolic in the other side with ${\rm Ker}F_p^2\cap {\rm Im}F_p^2=\{0\}$ and non-effectively hyperbolic on $S$ with ${\rm Ker}F_p^2\cap {\rm Im}F_p^2\neq \{0\}$.
\end{description}
\end{lem}
Let us denote
\[
\Sigma^{\pm}=\{(x,\xi)\in \Sigma\mid \pm(\{\xi_0,\phi_2\}^2-\{\phi_1,\phi_2\}^2)>0\}.
\]
Since the eigenvalues of $F_p$ are  $0$ and $\pm \sqrt{\{\xi_0,\phi_2\}^2-\{\phi_1,\phi_2\}^2}$ on $\Sigma$ so that $F_p$ has non-zero real eigenvalues on $\Sigma^{+}$ and non-zero pure imaginary eigenvalues on $\Sigma^{-}$ in the case (iii). Let us set
\[
2\kappa(\rho)^2=|\{\xi_0,\phi_2\}^2-\{\phi_1,\phi_2\}^2|
\]
and we make precise the meaning ``simply" in \eqref{eq:assump}, namely we assume that there is  $C>0$ such that
\begin{equation}
\label{eq:distance}
{\rm dist}_{\Sigma}(\rho,S)/C\leq \kappa(\rho)\leq C\,{\rm dist}_{\Sigma}(\rho,S)
\end{equation}
on $\Sigma$ where ${\rm dist}_{\Sigma}(\rho,S)$ denotes the distance from $\rho$ to $S$  on $\Sigma$. Our aim in this paper is to complete the proof of the following result:
\begin{thm}
\label{thm:main}
Assume \eqref{eq:assump} and that there is no bicharacteristic tangent to $\Sigma$ and there exist $\epsilon>0$, $C>0$ such that
\begin{equation}
\label{eq:Melcon}
(1-\epsilon)\mu(\rho)+ {\mathsf{Re}}P_{sub}(\rho)\geq \epsilon,\;\;|{\mathsf{Im}}P_{sub}(\rho)|\leq Ce(\rho),\;\rho\in \Sigma \cap\{|\xi|=1\}
\end{equation}
where $\pm e(\rho)$ $(e(\rho)\geq 0)$ are real eigenvalues and $\pm i\mu(\rho)$ $(\mu(\rho)\geq 0)$ are pure imaginary eigenvalues of $F_p(\rho)$. We also assume $n=2$ in the caae {\rm (iii)}. Then the Cauchy problem for $P$ is well-posed in any Gevrey class $\gamma^{(s)}$ for  $s> 1$. 
\end{thm}
The case (i) in Theorem \ref{thm:main}, namely $e(\rho)\equiv 0$ on $\Sigma$ was proved in \cite{BPP} while in \cite{Ni3}, it was proved under less restrictive assumption, the non existence of  bicharacteristics tangent to $S$. The case (ii) in Theorem \ref{thm:main} and hence $\mu(\rho)\equiv 0$ on $\Sigma$,   was proved in \cite{Ni4}. Some transition cases from effectively hyperbolic to non-effectively hyperbolic are studied in \cite{BB, BE, Es}. In particular in \cite{BE, Es} a typical case of (iii) was studied but the condition \eqref{eq:Melcon} was not investigated. In this paper we give a proof of Theorem \ref{thm:main} for the case (iii) assuming $n=2$, while if $n=1$ the case ${\rm Ker}F_p^2\cap {\rm Im}F_p^2\neq \{0\}$ never occur. 

\begin{rem}
\label{re:Mel} \rm For differential operators, the condition \eqref{eq:Melcon} with $\epsilon=0$ can be expressed as 
\[
{\rm dist}_{\C}\big(P_{sub}(\rho),[-{\rm Tr}^{+}F_p(\rho),{\rm Tr}^{+}F_p(\rho)]\big)\leq Ce(\rho)
\]
which generalizes the Ivrii-Petkov-H\"ormander condition (\cite{IP, Ho1}) and R.Melrose conjectured in \cite{Mel} that this condition is necessary for the Cauchy problem to be $C^{\infty}$ well-posed, but little is known about necessary conditions for the well-posedness when the spectral structure of $F_p$ changes. 
\end{rem}
\begin{rem}
\label{rem:bicha:e} \rm
With $X^{\pm}=\{\xi_0,\phi_2\}H_{\xi_0}-\{\phi_1,\phi_2\}H_{\phi_1}\pm \sqrt{2}\kappa(\rho) H_{\phi_2}$ it is easy to see 
\[
F_p(\rho)X^{\pm}=\pm e(\rho)X^{\pm},\quad \rho\in\Sigma^{+}
\]
and there exist exactly two bicharacteristics passing $\rho$ transversally to $\Sigma^{+}$ with tangents $X^{\pm}$ (see \cite{KoNi}).  Since $d\phi_2(X^{\pm})=\{\xi_0,\phi_2\}^2-\{\phi_1,\phi_2\}^2=2\kappa(\rho)^2>0$ this implies that the surface $\phi_2=0$ is spacelike on $\Sigma^{+}$. On the other hand there is no bicharacteristic reaching $\Sigma^{-}$(see \cite{Iv1}).
\end{rem}
%


\section{Idea of the proof of Theorem \ref{thm:main}}

 From Lemma \ref{lem:FfF} we have ${\rm Ker}F_p^2\cap {\rm Im}F_p^2\neq \{0\}$ on $\Sigma^{-}$ and there is no bicharacteristic tangent to $\Sigma^{-}$ by assumption. Then thanks to \cite[Theorem 3.3]{Ni1} $p$ admits an elementary decomposition microlocally at every point on $\Sigma^{-}$. As in \cite{BPP, Ni3} we try to decompose $p=-(\xi_0+\phi_1-\psi)(\xi_0-\phi_1+\psi)+q$ with $\psi=o(|\phi_1|)$ and non-negative $q$  verifying $\{\xi_0-\phi_1+\psi,q\}=O(q)$ in $\Sigma^{-}$. These requirements essentially determine $\psi$ and actually the non existence of tangent bicharacteristic assures that $\xi_0-\psi_1+\psi$ commutes against $q$ better than the usual case. On the other hand, as checked in Remark \ref{rem:bicha:e} the surface $\hphi_2=0$ is spacelike on $\Sigma^{+}$, then    \cite{Ni2, Ni4} suggests the use of  pseudodifferential weight $T\approx e^{\zeta\log{\hat \phi}_2}$ where $\zeta$ is a cutoff symbol to  $\Sigma^{+}$.  
 Our strategy for proving Theorem \ref{thm:main} is rather naive so that we make a such decomposition and  derive weighted energy estimates with the cutoff weight $T$. But the decomposition should be compatible with the cutoff weights and to achieve  this goal we must be careful in choosing cutoff symbols and in estimating errors caused by them.  The assumption $n=2$ enables us to choose all symbols which we need, including cutoff symbols, in $S_{3/4,1/2}$ and we carry out  pseudodifferential calculus within the framework of $S_{3/4,1/2}$ though we often need the calculus in smaller class than $S_{3/4,1/2}$.

In the rest of this section we express the assumptions in more explicit form. In what follows  we assume $n=2$ and we work in a conic neighborhood of ${\bar \rho}\in S$. Without restrictions we may assume ${\bar \rho}=(0,{\bf e}_3)$, ${\bf e}_3=(0,0,1)\in \R^{3}$ with a system of local coordinates  $x=(x_0,x')=(x_0,x_1,x_2)$. 
 From  (\ref{eq:distance}) and Lemma \ref{lem:FfF} one can write 
\begin{equation}
\label{eq:simple}
\{\xi_0,\phi_2\}^2-\{\phi_1,\phi_2\}^2=\te|\xi'|+c_1\phi_1+c_2\phi_2
\end{equation}
in a neighborhood of ${\bar\rho}$ where $S$ is defined by $\{\te=0\}\cap\Sigma$ and $d\te\neq 0$ on $S$ and hence $\Sigma^{\pm}=\Sigma\cap\{\pm \te>0\}$. Compare this to the case (i) and (ii) where we have $
\{\xi_0,\phi_2\}^2-\{\phi_1,\phi_2\}^2=\mp \te^2+c_1\phi_1+c_2\phi_2$ respectively (\cite{Ni3, Ni4}).  Here note that 
\[
e(\rho)=\left\{\begin{array}{cc}
\sqrt{2}\kappa(\rho)&\rho\in \Sigma^{+}\\
0&\rho\in \Sigma^{-}
\end{array}\right.,\quad 
\mu(\rho)=\left\{\begin{array}{cc}
0&\rho\in \Sigma^{+}\\
\sqrt{2}\kappa(\rho)&\rho\in \Sigma^{-}
\end{array}\right. .
\]
Since $\{\xi_0,\phi_2\}^2-\{\phi_1,\phi_2\}^2=\{\xi_0-\phi_1,\phi_2\}\{\xi_0+\phi_1,\phi_2\}=0$ on $S$ we may assume without restrictions that
\begin{equation}
\label{eq:ame}
\{\xi_0-\phi_1,\phi_2\}=0\quad\mbox{on}\quad S
\end{equation}
and $\{\xi_0,\phi_2\}=\{\phi_1,\phi_2\}>0$ on $S$ (see  \cite{Ni3, Ni4}). 
\begin{lem}
\label{lem:kame} In a conic neighborhood of ${\bar\rho}'=(0,{\bf e}_2)$ 
one can assume that
\[
\phi_2(x,\xi')=\hphi_2(x)e(x,\xi'),\;\;\te(x,\xi')|\xi_2|^{-1}=\psi(x')+f(x,\xi')\hphi_2(x)
\]
where $0\neq e(x,\xi')\in S^1_{1,0}$ and $f(x,\xi')\in S^0_{1,0}$. Moreover we have $\{\te,\phi_j\}=c_j\phi_2$ 
with $c_j\in S^0_{1,0}$.
\end{lem}
%
%
\noindent
Proof: Since $\{\xi_0,\phi_2\}\neq 0$ then one can write $\phi_2=(x_0-\psi_2(x,\xi'))b_2$ where $\psi_2$ is independent of $x_0$ and $b_2\neq 0$. From $\{\phi_1,\phi_2\}\neq 0$ we see $\{\psi_2,\phi_1\}\neq 0$.  This shows that $d\psi_2$ is not proportional to $\sum_{j=0}^2\xi_jdx_j$ at ${\bar\rho}$ because otherwise we would have $\phi_1(0,{\bf e}_2)=\dif \phi_1(0,{\bf e}_2)/\dif \xi_2\neq 0$. Since $\Xi_0=\xi_0$, $X_0=x_0$, $X_1=\psi_2$ verifies the commutation relations and $d\Xi_0$, $dX_0$, $dX_1$, $\sum_{j=0}^2\xi_jdx_j$ are linearly independent at ${\bar\rho}$, just observed above, these coordinates extends to homogeneous symplectic coordinates $(X,\Xi)$(see \cite[Theorem 21.1.9]{Ho1}). Switching the notation to $(x,\xi)$ we can assume that $\phi_2=(x_0-x_1)e$. Since $\{\phi_2,\phi_1\}\neq 0$ one can write $\phi_1=(\xi_1-\psi_1)b_1$ where $\psi_1$ is independent of $\xi_0$ and $\xi_1$. Writing $\psi_1(x,\xi_2)
={\bar \psi}_1(x',\xi_2)+e_1\phi_2$ and $\te(x,\xi')={\tilde \te}(x',\xi_2)+(x_0-x_1)\te_1+(\xi_1-{\bar \psi}_1)\te_2$ 
so that  $S$ is given by $\xi_0=0$, $x_0-x_1=0$, $\xi_1-{\bar\psi}_1(x',\xi_2)=0$, ${\tilde \theta}(x',\xi_2)=0$ 
where ${\tilde \te}(x',\xi_2)=\te(x_1,x_1,x_2,{\bar \psi}_1(x',\xi_2),\xi_2)$. Since ${\tilde \te}$ is of homogeneous of degree $1$ in $\xi_2$ one can write  
\[
{\tilde \theta}(x',\xi_2)={\tilde \theta}(x',1)\xi_2=\psi(x')\xi_2
\]
in a conic neighborhood of $(0,{\bf e}_2)$, {\it where we have used the assumption} $n=2$.
 Let us set $
{\te}=\psi(x')\xi_2+(\{\psi(x')\xi_2,\phi_1\}/\{\phi_1,\phi_2\})\phi_2$ 
then it is clear that $\{{\te} ,\phi_j\}=c_j\phi_2$ and hence this ${\te}$ is a desired one.
\qed
\begin{rem}
\label{rem:taku}\rm  Since the restriction $n=2$ is only used to prove Lemma \ref{lem:kame} then Theorem \ref{thm:main} is still true if we can choose a homogeneous symplectic coordinates such that Lemma \ref{lem:kame} holds.
\end{rem}
We now assume that $\phi_2$ and $\te$ satisfy Lemma \ref{lem:kame} and  set
\[
\hte=\te |\xi_2|^{-1},\quad \hphi_1=\phi_1|\xi'|^{-1}
\]
so that $\hte$ and ${\hat\phi}_1$ are homogeneous of degree $0$ in $\xi'$.
From (\ref{eq:ame}) we can write
\begin{equation}
\label{eq:abb}
\{\xi_0-\phi_1,{\hphi}_2\}={\hat c}\,\hte+c_1'{\hphi}_1+c_2'{\hphi}_2
\end{equation}
near ${ \bar\rho}$  where  ${\hat c}>0$ which follows from \eqref{eq:simple}. Since we have $\{\xi_0+\phi_1,\phi_2\}|{\hat c}\,\hte||e|=2\kappa^2$ on $\Sigma$ and $
\{\xi_0+\phi_1,\phi_2\}/2\{\phi_1,\phi_2\}=1$ 
on $S$ then for any $\ep>0$ there is a neighborhood of ${\bar\rho}$ where we have
\begin{equation}
\label{eq:abc}
(1-\ep)\kappa^2(\rho)\leq \{\phi_1,\phi_2\}|{\hat c}\,{\hte}||e|\leq (1+\ep)
\kappa^2(\rho).
\end{equation}
Here we examine how the non existence of tangent bicharacteristics reflects on the Poisson brackets of symbols .
\newtheorem{pro}{Proposition}[section]
\begin{pro}
\label{pro:nibicha} {\rm (\cite[Proposition 2.1]{Ni4})}
Assume $\{\theta,\phi_j\}=O(|\phi|)$ and that there is no bicharacteristic tangent to $\Sigma$. Then we have
\[
\{\xi_0,\te\}(\rho)= 0,\quad \{\{\xi_0-\phi_1,\phi_2\},\phi_2\}(\rho)= 0,\;\;\rho\in S.
\]
\end{pro}
\begin{lem}
\label{lem:a:3}
Assume that $\{\{\xi_0-\phi_1,\phi_2\},\phi_2\}=0$ on $S$. Then one can write $
\{\xi_0-\phi_1,\hphi_2\}={\hat c}\,\hte+c_0\hte \hphi_1+c_1\hphi_1^2+c_2\hphi_2$.
\end{lem}
\begin{lem}
\label{lem:a:4}
Assume that $\{\xi_0,\hte\}=0$, $\{\{\xi_0-\phi_1,\phi_2\},\phi_2\}=0$ on $S$. Then we have $
\{\xi_0-\phi_1,\hte\}=c_0\hte +c_1\hphi_1^2+c_2\hphi_2$.
\end{lem}
%
\noindent
Proof:
 Note that $\{\xi_0-\phi_1,\hte\}=\al \hte+\be\hphi_1+\ga\hphi_2$. On the other hand we see
\[
\{\hte,\{\xi_0-\phi_1,\phi_2\}\}=O(|\hphi|),\;\;\{\xi_0-\phi_1,\{\hte,\phi_2\}\}=O(|(\hte,\hphi)|).
\]
Then from the Jacobi identity it follows that $\be=O(|(\hte,\hphi)|)$ and hence we have $
\{\xi_0-\phi_1,\hte\}=\al \hte+c_0\hte\hphi_1+c_1\hphi_1^2+c_2\hphi_2$ 
which proves the assertion.
\qed
\begin{cor}
\label{cor:ame}
We have $\{\xi_0,\hte\}=c_0\hte+c_1\hphi_1^2+c_2\hphi_2$.
\end{cor}
%

\section{Cutoff and weight symbols}

We use the same notation as in \cite{Ni4}.
We first make a dilation of the coordinate $x_0$; $x_0\to \mu x_0$ with small $\mu>0$ so that $P(x,\xi,\mu)=\mu^2P(\mu x_0,x',\mu^{-1}\xi_0,\xi')$ will be
\begin{align*}
p(\mu x_0,x',\xi_0,\mu\xi')+\mu P_1(\mu x_0,x',\xi_0,\xi')+\mu^2P_0(\mu x_0,x')\\
=p(x,\xi,\mu)+P_1(x,\xi,\mu)+P_0(x,\mu).
\end{align*}
In what follows we often express such symbols dropping $\mu$. It is easy to see that $a(\mu x_0,x',\mu \xi')=a(x,\xi',\mu)\in S(\mxi^m,g_0)$ if $a(x,\xi')\in S^m_{1,0}$ where $g_0=|dx|^2+\xim^{-2}|d\xi'|^2$. To prove the well-posedness of the Cauchy problem, applying \cite[Theorem 1.1]{Ni5}, it suffices to derive energy estimates for $P_{\xi'}$ which coincides with original $P$ in a conic neighborhood of $(0,0,\xi')$, $|\xi'|=1$. Thus 
we can assume that the following conditions are satisfied globally;
\begin{equation}
\label{eq:global:assump}
\left\{\begin{array}{lll}
p(x,\xi)=-\xi_0^2+\phi_1(x,\xi')^2+\phi_2(x,\xi')^2,\;\;\phi_j\in S(\mxi,g_0),\\[3pt]
\{\xi_0,\phi_1\}=d_1\phi_1+d_2\phi_2,\;\;d_j\in \mu S(1,g_0),\\[3pt]
\{\xi_0-\phi_1,\hphi_2\}=\mu {\hat c}\,\hte+c_0\hte\hphi_1+c_1\hphi_1^2+c_2\hphi_2,\\[4pt]
\{\xi_0,\hte\}=c'_0\hte+c'_1\hphi_1^2+c'_2\hphi_2,\\[4pt]
\{\phi_1,\hphi_2\}\geq c\mu,\;\;c>0
\end{array}\right.
\end{equation}
where $c_j$, $c'_j$$\in \mu S(1,g_0)$ and  $\hte\in S(1,g_0)$ verifies
\begin{equation}
\label{eq:global:assump:b}
\{\hte,\phi_j\}=c_j\hphi_2,\;\;c_j\in \mu S(1,g_0)
\end{equation}
and 
 $\sup|\hte|$, $\sup|\hphi_j|$ can be assumed to be sufficiently small, shrinking a conic neighborhood of $(0,0,\xi')$ where we are working.
 
Let us put $P_{sub}=P^s_1+iP^s_2$ with real $P^s_i\in \mu S(\mxi,g_0)$ then from  \eqref{eq:Melcon} and \eqref{eq:abc} the following conditions can be assumed to be satisfied globally; 
\begin{equation}
\label{eq:GIPH}
\left\{\begin{array}{ll}
\mu^{1/2}\sqrt{{\hat c}\,\{\phi_1,\hphi_2\}|\hte|}|e|+P^s_1\geq c\mu\mxi,\; \hte<0,\; \;P^s_1\geq c\mu\mxi,\; \hte>0,\\[4pt]
P^s_2= \mu c_0\hte \mxi+c_{11}\phi_1+c_{12}\phi_2
\end{array}\right.
\end{equation}
with a constant $c>0$ and
 $c_0\in S(1,g_0)$,  $c_{ij}\in\mu S(1,g_0)$ where $c_0=0$ for $\hte<0$. Recall from \cite{Ni4}
\[
\left\{\begin{array}{lll}
\phi=\xim^{1/2}\big(\hphi_2+w\big),\\
\Phi=\pi+i\big\{\log{(\hphi_2+i\om)}-\log{(\hphi_2-i\om)}\big\}=\pi-2\arg{(\hphi_2+i\om)},\\
w=(\hphi_2^2+\xim^{-1})^{1/2},\quad
\om=(\hphi_1^4+\xim^{-1})^{1/2},\\
\rho^2=\hphi_2^2+\om^2=\hphi_2^2+\hphi_1^4+\xim^{-1}\geq (w^2+\omega^2)/2
\end{array}\right.
\]
where $\phi$ plays a major role in our arguments and $\Phi $ is introduced in order to manage the energy estimates in the region $C\hphi_1^2\geq w$. Note that
\begin{equation}
\label{eq:Fai}
\{F,\Phi\}=2(\omega\{F,\hphi_2\}-\hphi_2\{F,\omega\})/\rho^2.
\end{equation}
 We use the following   metrics
\[
\left\{\begin{array}{ll}g=w^{-2}|dx|^2+w^{-1}\xim^{-2}|d\xi'|^2,\\
g_1=(\rho^{-1}+\om^{-1/2})^2|dx|^2+\om^{-1}\xim^{-2}|d\xi'|^2,\\
{\tilde g}=(w^{-1}+\omega^{-1/2})^2|dx|^2+\xim^{-3/2}|d\xi'|^2,\\
\bg=\xim^{-1}|dx|^2+\xim^{-3/2}|d\xi'|^2.
\end{array}\right.
\]
Note that $g$, $g_1\leq {\tilde g}\leq \bg$ and $\bg$ is the metric defining the class $S_{3/4,1/2}$ for any fixed $\mu>0$. As checked in \cite{Ni4}, we have $\omega\in S(\omega,g_1)$, $\rho\in S(\rho,g_1)$ and $\Phi\in S(1,g_1)$. With cutoff symbol $\zeta(x,\xi')=\zeta(\hte w^{-1})$ we define the following  weight 
\begin{equation}
\label{eq:sukaju}
T=\exp{\big(n \zeta^2(\chi^2\log{\phi}+\Phi)\big)}
\end{equation}
where $\chi=\chi(\hphi_1^2w^{-1})$ and $\zeta(s)=1$ in $s\geq -b_1$ and $\zeta(s)=0$ in $s\leq -b_2$ with $\zeta'(s)\geq 0$ and $n$ is a positive parameter.


\begin{center}

\begin{tikzpicture}[scale=0.7, transform shape]


\draw (-9,0)--(-3,0);

\draw (-2,0)--(6,0);

\node at (-9.5,1) {\Large$1$};

\draw[rounded corners=8pt] (-9,1)--(-7,1)--(-6.3,0)--(-3,0);
\draw[rounded corners=8pt] (-9,0)--(-5.7,0)--(-5,1)--(-3,1);

\node at (-6,-0.4) {\Large$d_2$};

\node at (-5,-0.4) {\Large$d_3$};

\node at (-4,1.4) {\Large$\chi_2$};

\draw[dashed] (-4.9,0)--(-4.9,1);

\draw[dashed] (-8,0)--(-8,1);

\node at (-7,-0.4) {\Large$d_1$};

\node at (-8,-0.4) {\Large$0$};

\node at (-7.5,1.4) {\Large$\chi$};

\draw[dashed] (-7.1,0)--(-7.1,1);

\node at (6.5,1) {\Large$1$};

\draw[rounded corners=8pt] (-2,1)--(-1,1)--(-0.3,0)--(6,0);

\draw[rounded corners=8pt] (-2,0)--(0.3,0)--(1,1)--(6,1);

\draw[rounded corners=8pt] (-2,0)--(4,0)--(5,1)--(6,1);

\node at (-1.2,-0.4) {\Large$-b_3$};

\node at (-1.2,1.4) {\Large${\zeta}_{-}$};

\draw[dashed] (-1.2,0)--(-1.2,1);

\node at (-0.2,-0.4) {\Large$-b_2$};

\node at (4,0.5) {\Large$\zeta_{+}$};

\node at (0.9,-0.4) {\Large$-b_1$};

\draw[dashed] (1.1,0)--(1.1,1);

\node at (2,-0.4) {\Large$0$};

\node at (2.5,1.4) {\Large$\zeta$};

\draw[dashed] (2,0)--(2,1);

\node at (3,-0.4) {\Large$b_1$};


\node at (4,-0.4) {\Large$b_2$};

\node at (5,-0.4) {\Large$b_3$};

\draw[dashed] (5.1,0)--(5.1,1);

\end{tikzpicture}
\end{center}

Let $\zeta_{\pm}(x,\xi')=\zeta_{\pm}(\hte w^{-1})$ and $\chi_2(x,\xi')=\chi_2(\hphi_1^2w^{-1})$ where  $\zeta_{\pm}(s)=1$ in $\pm s\geq b_3$ and $0$ in $\pm s\leq b_2$ so that $\zeta\zeta_{+}=\zeta_{+}$ and $\zeta\zeta_{-}=0$. We simply write $\chi$, $\chi_2$ for $\chi(x,\xi')$ and $\chi_2(x,\xi')$ and $\zeta$, $\zeta_{\pm}$ for $\zeta(x,\xi')$ and $\zeta_{\pm}(x,\xi')$ if there is no confusions.  It is easy to check $\chi$, $\chi_2\in S(1,g)$. 
  As for new cutoff symbols $\zeta$, $\zeta_{\pm}$ we have
\begin{lem}
\label{lem:ameagari} Let $G=w^{-2}|dx|^2+\xim^{-2}|d\xi'|^2\;(\leq g)$ then  
$w\in S(w,G)$ and $\phi\in S(\phi,G)$. We have also 
$\zeta,\zeta_{\pm}\in S(1,G)$.  Let $s\in\R$ then $\zeta_{+}\hte^s \in S(|\hte|^s,G)$. Moreover if $0<s\leq 1$ and $|\al|\neq 0$ we have $|(\zeta_{+}\hte^s)^{(\al)}_{\be)}|\leq C_{\al\be}w^s\xim^{-|\al|}w^{-|\be|}$. 
\end{lem}
%
\noindent
Proof:
 To prove $\phi\in S(\phi,G)$, with ${\tilde  \phi}=\hphi_2+w$, it is enough to show ${\tilde \phi}\in S({\tilde \phi},G)$. Note that one can write
\[
\dif_x^{\be}\dif_{\xi'}^{\al}{\tilde\phi}=\frac{\dif_x^{\be}\dif_{\xi'}^{\al}\hphi_2(x)}{w}{\tilde\phi}+\frac{\dif_x^{\be}\dif_{\xi'}^{\al}\xim^{-1}}{2w}=b_{\al\be}{\tilde \phi}+a_{\al\be}
\]
with $b_{\al\be}\in S(w^{-|\be|}\xim^{-|\al|},G)$ and $a_{\al\be}\in S((w^{-1}\xim^{-1})w^{-|\be|}\xim^{-|\al|},G)$ for $|\al+\be|=1$. By induction on $|\al+\be|$ we see easily $
\dif_x^{\be}\dif_{\xi'}^{\al}{\tilde\phi}=b_{\al\be}{\tilde \phi}+a_{\al\be}$ 
with $b_{\al\be}\in S(w^{-|\be|}\xim^{-|\al|},G)$ and $a_{\al\be}\in S((w^{-1}\xim^{-1})w^{-|\be|}\xim^{-|\al|},G)$ for any $\al$, $\be$. Since 
$w^{-1}\xim^{-1}\leq 2{\tilde\phi}$ we get the assertion.
To prove $\zeta\in S(1,G)$ it suffices to show that
\begin{equation}
\label{eq:kaisi}
|\zeta'\dif_x^{\be}\dif_{\xi'}^{\al}(\hte w^{-1})|\leq C_{\al\be}w^{-|\be|}\xim^{-|\al|}.
\end{equation}
By Lemma \ref{lem:kame} without restrictions we may assume $
\hte(x,\xi')=\psi(x')+f(x,\xi')\hphi_2(x)$ from which it follows $|\dif_{\xi'}^{\al}\hte|\leq C_{\al}\xim^{-|\al|}w$ for $|\al|\geq 1$. Noting  $|\zeta'\hte w^{-1}|\leq C$ we get \eqref{eq:kaisi}. On the support of $\zeta_{+}$ the estimate 
\[
|(\hte^s)^{(\al)}_{(\be)}|\leq \sum C_{\al_1,\ldots,\be_k}\hte^s|\hte^{(\al_1)}_{(\be_1)}|\hte^{-1}\cdots|\hte^{(\al_k)}_{(\be_k)}|\hte^{-1}
\]
holds where $|\al_i+\be_i|\geq 1$ and $\al_1\cdots+\al_k=\al$, $\be_1\cdots+\be_k=\be$. On the other hand  Lemma \ref{lem:kame} shows that $|\hte^{(\al_i)}_{(\be_i)}|\leq C_{\al_i\be_i}\xim^{-|\al_i|}w^{1-|\be_i|}$ if $|\al_i|\neq 0$ and bounded by $ C_{\be_i}$ if $|\al_i|=0$. Since $\hte^{-1}w$ is bounded on the support of $\zeta_{+}$  the third assertion is clear. If $|\al_i|\neq 0$ then noting $\hte^s|\hte^{(\al_i)}_{(\be_i)}\hte^{-1}|\leq C_{\al_i\be_i}w^s\xim^{-|\al_i|}w^{-|\be_i|}$ on the support of $\zeta_{+}$ one gets the last assertion.
\qed
\begin{rem}\rm  If $n>2$ the $\psi(x')$ in Lemma \ref{lem:kame} would depend on $\xi'$ also and hence $\zeta$, $\zeta_{\pm}$ does not belong to $S(1,g)$ in general. 
\end{rem}

To decompose $p$ let us define
\begin{equation}
\label{eq:pusai}
\psi=(-{ h}\zeta^2_{-}+\nu\zeta^2_{+}) \hte\phi_1+\chi_2\phi_1^3\mxi^{-2}={\tilde \zeta}\hte\phi_1+\chi_2\phi_1^3\mxi^{-2}
\end{equation}
with a  positive parameter $0<\nu\ll 1$ which   will be determined later where ${\tilde\zeta}=-{ h}\zeta^2_{-}+\nu\zeta^2_{+}$ with ${ h}=\mu{\hat c}\{\phi_1,\hphi_2\}^{-1}>0$. Using $\psi$ we rewrite $p$ as
\begin{equation}
\label{eq:ba}
\begin{split}
&p=-(\xi_0+\phi_1-\psi)(\xi_0-\phi_1+\psi)+2\psi\phi_1-\psi^2+\phi_2^2\\
&=-(\xi_0+\phi_1-\psi)(\xi_0-\phi_1+\psi)+q
\end{split}
\end{equation}
where 
 \[
 \left\{\begin{array}{ll}
 q=\phi_2^2+2a^2{\tilde\zeta}  \hte\phi_1^2+2a^2\chi_2 \phi_1^4\mxi^{-2},\\[3pt]
 a=(1-{\tilde\zeta}\hte/2-\chi_2\phi_1^2\mxi^{-2}/2)^{1/2}.
\end{array}\right.
\]
The main part of $\{\xi_0-\phi_1+\psi,q\}$ will be $\{\xi_0-\phi_1+\psi,\phi_2^2\}$ which is required to be $O(q)$ in $\te<0$ as explained above. Indeed, by our choice, we have
\begin{equation}
\label{eq:bunkai}
\begin{split}
\{\xi_0-\phi_1+\psi,\hphi_2\}
=\mu(1-\zeta_{-}^2){\hat c}\,\hte+\mu\nu {\hat h}^{-1}{\hat c}\,\zeta^2_{+}\hte\\+c_1\hphi_1^2+c_2\hte\hphi_1+c_3\hphi_2
\end{split}
\end{equation}
where $1-\zeta_{-}^2=0$ in $\hte\leq -b_3w$ so that $|(1-\zeta_{-}^2)\hte|\leq Cw$ in $\hte\leq 0$.
\begin{lem}
\label{lem:toyama}
We have $({\tilde \zeta}\hte)^{(\al)}_{(\be)}$, $(\chi_2\hphi_1^2)^{(\al)}_{(\be)}\in S(\xim^{-|\al|}, g)$ for $|\al+\be|=1$. Hence the same holds for $a^{(\al)}_{(\be)}$.  In particular $|({\tilde \zeta}\hte)^{(\al)}_{(\be)}|$, $|(\chi_2\hphi_1^2)^{(\al)}_{(\be)}|$ and $|a^{(\al)}_{(\be)}|$ are bounded by $C_{\al\be}w^{1/2}\xim^{-|\al|}w^{-|\al|/2-|\be|}$ for $|\al+\be|\geq 1$ and bounded by $C_{\al\be}w\xim^{-|\al|}w^{-|\al|/2-|\be|}$ for $|\al+\be|\geq 2$.
\end{lem}
In this paper $Op(\phi)$ denotes the Weyl quantized pseudodifferential operator with symbol $\phi$ and we denote $Op(\phi)Op(\psi)=Op(\phi\#\psi)$. We often use the same letter to denote a symbol and the operator with such symbol if there is no confusion. Thus we denote
\[
Op(\phi\psi)u=\phi\psi u,\quad Op(\phi)Op(\psi)u=\phi(\psi u).
\]
We make some additional preparations (see \cite{Iv2}). Let $c=id_1+ic_{11}$ with $c_{11}$, $d_1$ in \eqref{eq:GIPH},  \eqref{eq:global:assump:b} and we set $M=\xi_0+\phi_1-\psi+c$, $\La=\xi_0-\phi_1+\psi-c$ and write 
\[
p+P_1^s+iP_2^s=-M\#\Lambda+Q=-M\#\Lambda+q+T_1+iT_2.
\]
Note that $-(\xi_0+\phi_1-\psi)(\xi_0-\phi_1+\psi)=
-M\Lambda-c\phi_1-2c \psi-c^2$.
 In view of Lemma \ref{lem:toyama}  it is not difficult to check 
\[
M\#\Lambda=M\Lambda+i\{\xi_0,\phi_1-\psi+c_{11}\}+c_1w^{1/2}\phi_1+c_2\hphi_1^2\mxi+R
\]
 with $c_i\in \mu S(1,\bg)$, $R\in \mu^2 S(w^{-1},\bg)$. Therefore we see from \eqref{eq:global:assump} that $T_1$ satisfies 
\begin{equation}
\label{eq:mao}
\left\{\begin{array}{ll}
 \mu^{1/2}\sqrt{{\hat c}\,\{\phi_1,\hphi_2\}|\hte|}|e(x,\xi')|+T_1\geq 2{\bar\kappa}\mu\mxi,\quad \hte<0,\\[4pt]
T_1\geq 2{\bar\kappa}\mu\mxi,\quad \hte>0
\end{array}\right.
\end{equation}
with some ${\bar\kappa}>0$ and $T_2$ can be written
\begin{equation}
\label{eq:ToneT}
T_2=\mu  c_0\hte\mxi+b_0\hte\phi_1+b_1\hphi_1^2\mxi+b_2\phi_2+b_3w^{1/2}\phi_1
\end{equation}
with $b_i\in S(1,\bg)$. Thus $T_2\mxi^{-1}=O(|(\hte,\hphi_1^2,\hphi_2, w^{1/2}\hphi_1)|)$, a linear combination without  $\hphi_1$. We transform $P$ by $T$ so that 
\[
PT=T{\tilde P},\quad {\tilde P}=-{\tilde M}{\tilde \La}+{\tilde Q}.
\]
To simplify notations we set $\Psi=\zeta^2(\chi^2\log{\phi}+\Phi)$. Then we have
\begin{lem}
\label{lem:TT}
We have $
T=e^{n\Psi}\in S(e^{n \Psi},(\log^2{\xim})\bg)$.
\end{lem}
%
\noindent
Proof:
 Note that $\dif_{x}^{\be}\dif_{\xi'}^{\al}\log{\phi}=\phi^{-1}\dif_x^{\be}\dif_{\xi'}^{\al}\phi$ and 
$\phi^{-1}\in S(\phi^{-1},g)$ for $|\al+\be|=1$. Since $|\log{\phi}|\leq C\log{\xim}$ and $g$, $g_1\leq \bg$ the assertion is clear. 
\qed
 
Let us write ${\tilde M}=D_0-{\tilde m}(x,D')$, ${\tilde \La}=D_0-{\tilde\la}(x,D')$ and fix any small $\varep>0$.
\begin{pro}{\rm (\cite{Ni3, Ni4})}
\label{pro:ener:id} Let ${\tilde P}=-({\tilde M}-i\gamma\la_{\mu}^{2\epsilon})({\tilde \Lambda}-i\gamma\la_{\mu}^{2\epsilon})+{\tilde Q}$ then we have
\begin{equation}
\label{eq:minore}
\begin{split}
2{\mathsf{Im}}({\tilde P}u,{\tilde\La}u)\geq \frac{d}{dx_0}(\|{\tilde\La}u\|^2+(({\mathsf{Re}}\,{\tilde Q})u,u)+\ga^2\|\Dm^{2\varep}u\|^2)\\
+\ga\|\la_{\mu}^{\varep}({\tilde\La}u)\|^2
+2\ga {\mathsf{Re}}(\la_{\mu}^{2\varep}({\tilde Q}u),u)+2(({\mathsf{Im}}\,{\tilde m}){\tilde\La}u,{\tilde\La}u)\\
+2{\mathsf{Re}}({\tilde\La}u,({\mathsf{Im}}\,{\tilde Q})u)
+{\mathsf{Im}}([D_0-{\mathsf{Re}}{\tilde \la},{\mathsf{Re}}\,{\tilde Q}]u,u)
\\
+2{\mathsf{Re}}(({\mathsf{Re}}\,{\tilde Q})u,({\mathsf{Im}}{\tilde \la} )u)
+\frac{\ga^3}{2}\|\la_{\mu}^{3\varep}u\|^2
+2\ga^2 (\la_{\mu}^{4\varep}({\mathsf{Im}}{\tilde \la}) u,u).
\end{split}
\end{equation}
\end{pro}
In this paper positive large parameters $n$, $\gamma$ and a positive  small parameter $\mu$ are assumed to satisfy $
n\mu^{1/4}\ll 1$ and $ \gamma \mu^{4}\gg 1$.
\begin{rem}\rm
The weight $\lr{\mu D'}^{2\varep}$ is introduced to control error terms $\log^N{\lr{ D'}}$, caused by metric $(\log^2{\xim})\bg$, and hence we can choose $\varep>0$ as small as we please, which determines the well-posed Gevrey class $\gamma^{(1/2\varep)}$. Actually the Cauchy problem is well-posed in the space consisting of  all $C_0^{\infty}$ functions with Fourier transform bounded by $\exp{\big(-C\log^N\lr{\xi'}\big)}$ 
with some $C>0$, $N>0$.
\end{rem}

\begin{defin}
\label{dfn:pluszero}\rm  We set $\lambda=\mxi$, $\lambda_{\mu}=\xim$ and $\la^{s+0}=\mxi^s\xim^{+0}$. 
We write $a\in S(\la_{\mu}^{s+0},g)$ ($a\in S(\la^{s+0},g)$) 
if $a\in S(\xim^{s+\varep},g)$ ($a\in S(\mxi^s\xim^{\varep},g)$ for any $\varep>0$. We also denote
\[
\|Au\|\leq C\|\la_{\mu}^{s+0}u\|\;\;(\|\la^{s+0}u\|)
\]
if $\|A u\|\leq C\|\Dm^{s+\varep}u\|$ ($\|A u\|\leq C\|\mD^s\Dm^{\varep}u\|$) for any $\varep>0$ with some $C>0$ independent of $\varep>0$.
\end{defin}
%
%

%

\section{Transformed symbols ${\tilde \la}$, ${\tilde m}$}

We first list up several properties of cutoff symbols. 
\begin{lem}
\label{lem:innsi}
We have
\begin{equation}
\label{eq:katta}
\begin{split}
\chi\chi_2=0,\;\;\zeta\zeta_{-}=0,\;\;\zeta\zeta_{+}=\zeta_{+},\;\;{\tilde\zeta}\zeta=\nu\zeta_{+}^2,\\
\hphi_2,\;\;\chi \hphi_1^2,\;\;\zeta' \hte,\;\;\zeta_{\pm}'\hte \in S(w,g),\;\;\chi\hphi_1\in S(w^{1/2},g),\\
(1-\zeta^2_{-}-\zeta^2_{+})\hte,\;\;\zeta (1-\zeta^2_{+})\hte\in S(w,g)
\end{split}
\end{equation}
where $\zeta'=\zeta'(\hte w^{-1})$. We also have $
\{\chi,\la_{\mu}^s\}$, $\{\zeta,\la_{\mu}^s\}\in S(w^{-1}\la_{\mu}^{s-1},g)$.
\end{lem}

Denote $W^{\al}_{\be}=T^{-1}\dif_x^{\be}\dif_{\xi'}^{\al}T$ and note that we have for $a\in S(\la_{\mu}^{s+0} w^t,g)$ or  $a\in S(\la^s,g_0)$
\begin{eqnarray*}
&&a\#T=T\#a-i nT\{a,\Psi\}\\
&&+\frac{i}{8}\;T\!\sum_{|\al+\be|=3}\frac{(-1)^{|\be|}}{\al!\be!}\big(a^{(\al)}_{(\be)}W^{\be}_{\al}-W^{\al}_{\be}a^{(\be)}_{(\al)}\big)+T\#R
\end{eqnarray*}
 with some $R\in  S(w^t\la_{\mu}^{s-5/4+0},\bg)$ or $R\in  S(\la^s\la_{\mu}^{-5/2+0},\bg)$ respectively. From Lemma \ref{lem:toyama} it follows that $\psi^{(\al)}_{(\be)}W^{\be}_{\al}\in S(1,\bg)$ for $|\al+\be|=3$ then the main parts of ${\mathsf{Im}}\,{\tilde m}$ and ${\mathsf{Im}}\,{\tilde \la}$ are, up to the parameter $n$ 
\begin{align*}
\{\xi_0\pm \phi_1\mp \psi,\Psi\}
=\zeta^2\{\xi_0\pm\phi_1\mp\psi,\chi^2\log{\phi}+\Phi\}\\
+\{\xi_0\pm \phi_1\mp \psi,\zeta^2\}(\chi^2\log{\phi}+\Phi).
\end{align*}
To estimate $\{\xi_0\pm\phi_1\mp\psi,\chi^2\log{\phi}+\Phi\}$ it suffices to repeat similar arguments as in \cite{Ni4} to get
\begin{equation}
\label{eq:nisiyama}
\{\xi_0\pm\phi_1\mp\psi,\chi^2\log{\phi}+\Phi\}=\{\xi_0\pm\phi_1\mp\psi,\hphi_2\}(r+2\omega \rho^{-2})+R
\end{equation}
with $R\in S(\la_{\mu}^{+0},\bg)$ where
\begin{align*}
&0\leq  r=\chi^2 w^{-1}+\de \in S(w^{-1}\la_{\mu}^{+0},g),\\
& 0\leq \delta=-2\chi\chi'\hphi_1^2w^{-3}\hphi_2\log{\phi}\in S(w^{-1}\la_{\mu}^{+0},g)
\end{align*}
and the fact $\delta\geq 0$ follows from \cite[Lemma 3.6]{Ni4} which was a key point to treat \eqref{eq:nisiyama}. We check how the term $\{\xi_0\pm \phi_1\mp \psi,\zeta^2\}(\chi^2\log{\phi}+\Phi)$ can be managed. It is not difficult to see
\begin{equation}
\label{eq:arata:b}
\begin{split}
& \{\xi_0\pm\phi_1\mp\psi,\zeta^2\}(\chi^2\log{\phi}+\Phi)\\
&=-2\zeta\zeta'\hte w^{-3}\hphi_2(\chi^2\log{\phi}+\Phi)\{\xi_0\pm\phi_1\mp\psi,\hphi_2\}+R
\end{split}
\end{equation}
with $R\in \mu S(\la_{\mu}^{+0},\bg)$.
Here we note
\begin{lem}
\label{lem:phipos}
We have 
\[
0\leq \Delta=-2\zeta\zeta'\hte w^{-3}\hphi_2(\chi^2\log{\phi}+ \Phi)\in S(w^{-1}\la_{\mu}^{+0},\bg).
\]
\end{lem}
%
\noindent
Proof:
Since $\hphi_2\log{\phi}\geq 0$ by \cite[Lemma 3.6]{Ni5} it is clear $
0\leq -2\chi^2\zeta\zeta'\hte w^{-3}\hphi_2\log{\phi}\in S(w^{-1}\la_{\mu}^{+0},g)$ because $\zeta'(\hte w^{-1})\hte\leq 0$. Noting that $0\leq \Phi=\pi-2\arg{(\hphi_2+i\om)}\leq \pi$ if $\hphi_2\geq 0$ and $-\pi\leq \Phi=\pi-2\arg{(\hphi_2+i\om)}\leq 0$ 
for $\hphi_2\leq 0$ it is also clear $\hphi_2\Phi\geq 0$ and hence $
0\leq -2\zeta\zeta'\hte w^{-3}\hphi_2 \Phi\in S(w^{-1},\bg)$. 
Thus we get the assertion.
\qed

 To simplify notations we set $\Gamma=r+2\omega\rho^{-2}$. From \eqref{eq:nisiyama} and \eqref{eq:arata:b} it suffices to consider $
n(\Delta+\zeta^2\Gamma)\{\xi_0\pm \phi_1\mp\psi,\hphi_2\}$. 
 As in \cite{Ni4} we set
\[
\left\{\begin{array}{ll}
e_1=\mu {\hat c}+\nu\{\phi_1,\hphi_2\},\;\;e_3=\{\xi_0+\phi_1,\hphi_2\},\\
e_2=\{\xi_0+\phi_1,\hphi_2\}-\nu \hte\{\phi_1,\hphi_2\}\zeta^2_{+}.
\end{array}\right.
\]
Noting Lemma \ref{lem:innsi} it is easy to see
\begin{equation}
\label{eq:oto}
\begin{split}
\{\xi_0-\phi_1+\psi,\hphi_2\}=\mu{\hat c}\hte+{\tilde \zeta}\{\phi_1,\hphi_2\}\hte+c_0\hte\hphi_1+3\chi_2\hphi_1^2\{\phi_1,\hphi_2\},\\
\{\xi_0+\phi_1-\psi,\hphi_2\}=\{\xi_0+\phi_1,\hphi_2\}-{\tilde \zeta}\{\phi_1,\hphi_2\}\hte-3\chi_2\hphi_1^2\{\phi_1,\hphi_2\}
\end{split}
\end{equation}
modulo $S(w,\bg)$. Noting $\zeta=\zeta_{+}^2+\zeta(1-\zeta_{+}^2)$, $\zeta(1-\zeta_{+}^2)\hte\in S(w,\bg)$ we see
\[
\left\{\begin{array}{ll}
\zeta^2\{\xi_0-\phi_1+\psi,\hphi_2\}=(e_1+a_1\hphi_1)\zeta^2_{+}\hte+a_2\zeta\hphi_1^2,\\[3pt]
\zeta^2\{\xi_0+\phi_1-\psi,\hphi_2\}=e_2\zeta^2+a_3\zeta\hphi_1^2
\end{array}\right.
\]
with $a_i\in S(1,\bg)$ modulo $S(w,\bg)$. Since $\Delta\hte$, $\Gamma\hphi_1^2\in S(\la_{\mu}^{+0},\bg)$ by Lemma \ref{lem:innsi}
 we see ${\mathsf{Im}}\,{\tilde\lambda}=n\zeta_{+}^2\Gamma(e_1+a\hphi_1)+R$ 
with $R\in S(\la_{\mu}^{+0},\bg)$ and $a\in S(\la_{\mu}^{+0},\bg)$. Similarly we have ${\mathsf{Im}}\,{\tilde m}=n\Delta (e_3+a'\hphi_1^2)+n\zeta^2\Gamma e_2+R$ with $R\in S(\la_{\mu}^{+0},\bg)$. Noting that the main part of ${\mathsf{Re}}\,{\tilde \lambda}$ comes from $\{\{\xi_0-\phi+\psi,\Psi\},\Psi\}$ we summarize
\begin{lem}
\label{lem:whatla}
We have
\[
\left\{\begin{array}{lll}
{\mathsf{Im}}\,{\tilde \la}=n (e_1+b_1\hphi_1)\Gamma\zeta^2_{+}\hte+R_1,\\
{\mathsf{Re}}\,{\tilde \la}=\phi_1-\psi+n (b_2\hte+b_3\hphi_1^2) w^{-1/2}+R_2,\\
{\mathsf{Im}}\,{\tilde m}=ne_2\zeta^2\Gamma+n(e_3+b_4\hphi_1^2)\Delta+R_3
\end{array}\right.
\]
where $b_i$, $R_i\in S(\la_{\mu}^{+0},\bg)$.
\end{lem}
\begin{lem}
\label{lem:LaM}
There  exists $c>0$ which is independent of $\nu>0$ such that we have
\begin{eqnarray*}
C({\mathsf{Im}}{\tilde \la}\,u,u)\geq c\mu n(\Gamma \zeta^2_{+}\hte u,u)
-C_1\|\la_{\mu}^{+0}u\|^2\\
\geq c\mu n(\Gamma(\zeta_{+}\hte^{1/2})u, (\zeta_{+}\hte^{1/2})u)-C_2\|\la_{\mu}^{+0}u\|^2,
\\[3pt]
C({\mathsf{Im}}\,{\tilde m}\,u,u)\geq c\mu n((\zeta^2\Gamma+\Delta)u,u)-C_4\|\la_{\mu}^{+0}u\|^2\\
\geq c\mu n(\Gamma(\zeta  u),\zeta u)+c\mu n(\Delta u,u)-C_5\|\la_{\mu}^{+0}u\|^2.
\end{eqnarray*}
We have also
\begin{eqnarray*}
C({\mathsf{Im}}{\tilde\la}\,u,u)
\geq c\mu n(\|\chi\zeta_{+}\hte^{1/2}w^{-1/2}u\|^2
+\|\zeta_{+}\hte^{1/2}\rho^{-1/2} u\|^2),\\[3pt]
C({\mathsf{Im}}\,{\tilde m}\,u,u)
\geq c\mu n(\|\zeta \chi w^{-1/2}u\|^2+\|\zeta\rho^{-1/2}u\|^2)
\end{eqnarray*}
modulo $C\|\la_{\mu}^{+0}u\|^2$ with some $C$, $C'>0$ independent of $\mu$.
\end{lem}
%
\noindent
Proof:
Since $\hphi_1(0,{\bf e}_2)=0$ we may assume ${\tilde e}_1=e_1+b_1\hphi_1\geq \mu c_1>0$.  Take $M>0$ so that $M{\tilde e}_1\geq \mu$.
Since $0\leq (M{\tilde e}_1-\mu)\Gamma\zeta^2_{+}\hte\in \mu  S(w^{-1}\la_{\mu}^{+0},\bg)\subset \mu S^{1/2+0}_{3/4,1/2}$ then from the Fefferman-Phong inequality (see \cite[Theorem 18.6.8]{Hobook}) it follows that
\begin{align*}
M({\tilde e}_1\Gamma\zeta^2_{+}\hte u,u)
\geq \mu (\Gamma \zeta^2_{+}\hte u,u)-C_1\|\la_{\mu}^{+0}u\|^2.
\end{align*}
Here note that $
\Gamma \zeta_{+}^2\hte=(\zeta_{+}\hte^{1/2})\#\Gamma\#(\zeta_{+}\hte^{1/2})+R$ 
with $R\in S(\la_{\mu}^{+0},\bg)$. 
Since $|(R u,u)|\leq C'\|\la_{\mu}^{+0} u\|^2$ for $R\in S(\la_{\mu}^{+0},\bg)$ the first assertion follows.  To show the second assertion it suffices to repeat the same arguments proving the first assertion.

To prove the third assertion we first note that 
\[
(\delta\zeta_{+}^2\hte u,u),\; (\zeta^2\delta u,u),\;(\Delta u,u)\geq -C\|\la_{\mu}^{+0}u\|^2
\]
which follows the Fefferman-Phong inequality since $\delta $, $\Delta\in S^{1/2+0}_{3/4,1/2}$ are non-negative. We then write $\chi^2\zeta_{+}^2\hte w^{-1}=\chi\zeta_{+}\hte^{1/2}w^{-1/2}\#\chi\zeta_{+}\hte^{1/2}w^{-1/2}+R$ with $R\in S(\la_{\mu}^{+0},\bg)$ because $\zeta_{+}\hte^{1/2}\in S(1, G)\subset S(1,\bg)$ by Lemma \ref{lem:ameagari} which gives the first term on the right-hand side. To get the second term on the right-hand side we note that on the support of $1-\chi^2$ we have $C\omega\geq \rho$ with some $C>0$ and  it is obvious that $w^{-1}\geq \rho^{-1}$. Therefore it follows $C(\chi^2\zeta_{+}^2\hte w^{-1}+\zeta_{+}^2\hte\omega\rho^{-2})\geq \zeta_{+}^2\rho^{-1}\hte$.  Then the Fefferman-Phong inequality gives
\[
C(\chi^2\zeta_{+}^2\hte w^{-1}u,u)+C(\zeta_{+}^2\hte \omega\rho^{-2}u,u)\geq \|\zeta_{+}\rho^{-1/2}\hte u\|^2-C\|\la_{\mu}^{+0}u\|^2
\]
which gives the second term. The proof of the last assertion is similar.
\qed
 
Applying Lemma \ref{lem:LaM} one can show
\begin{pro}
\label{pro:akasaka} We have
\begin{align*}
2(({\mathsf{Im}}\,{\tilde m}){\tilde\La}u,{\tilde\La}u)\geq 
c\mu n((\Gamma+\Delta)(\zeta\tLa u),(\zeta\tLa u))+c\mu n\|\chi\zeta w^{-1/2}\tLa u\|^2\\
+c\mu n\|\zeta \omega^{1/2}\rho^{-1}\tLa u\|^2+c\mu n\|\zeta\rho^{-1/2}{\tilde \Lambda}u\|^2-C\|\la_{\mu}^{+0}\tLa u\|^2
\end{align*}
with some $c>0$ independent of $\nu>0$ and some $C>0$. 
\end{pro}
%

\section{Estimate $\|{\tilde\Lambda}u\|$}

We first remark the following lemma which is easily checked using \eqref{eq:global:assump} and \eqref{eq:global:assump:b}.
\begin{lem}
\label{lem:DGB} Let ${\hat\zeta}$, ${\hat\chi}\in C^{\infty}(\R)$ such that ${\hat\zeta}'$, ${\hat\chi}'\in C_0^{\infty}(\R)$. Set ${\hat\zeta}={\hat\zeta}(\hte w^{-1})$ and ${\hat\chi}={\hat\chi}(\hphi_1^2w^{-1})$. Then 
we have
\begin{align*}
\{\hphi_1,{\hat\zeta}\}\in S(w^{-1}\la_{\mu}^{-1},\bg),\;\{\hphi_1,{\hat\chi}\}\in S(w^{-1}\la_{\mu}^{-1},\bg),\\
\{\hphi_1,w^{-1}\}\in S(w^{-2}\la_{\mu}^{-1},\bg),\;\{\hphi_1,\omega\rho^{-2}\}\in S(\rho^{-2}\la_{\mu}^{-1},\bg),\;\\
\{\hphi_2,{\hat\chi}\}\in S(w^{-1/2}\la_{\mu}^{-1},\bg),
\{\hphi_2,{\hat\zeta}\}\in S(\la_{\mu}^{-1},\bg),\\
\{\hphi_2,\omega\rho^{-2}\}\in S(\rho^{-3/2}\la_{\mu}^{-1},\bg),\;
\{\hphi_2,w^{-1}\}\in S(w^{-1}\la_{\mu}^{-1},\bg),\\
\{{\hat\zeta},\hte\}\in S(\la_{\mu}^{-1},\bg),\;\{{\hat\zeta},w^{-1}\}\in S(w^{-2}\la_{\mu}^{-1},\bg),\;\{{\hat\zeta},\omega\rho^{-2}\}\in S(\rho^{-3/2}w^{-1}\la_{\mu}^{-1},\bg),\\
\{{\hat\chi}, \hte\}\in S(\la_{\mu}^{-1},\bg),\;
\{{\hat\chi},w^{-1}\}\in S(w^{-1/2},\bg),\;\{{\hat\chi},\omega\rho^{-2}\}\in S(\rho^{-5/2}\la_{\mu}^{-1},\bg),\\
\{{\hat\zeta},{\hat\chi}\}\in S(w^{-3/2}\la_{\mu}^{-1},\bg),\;
\{w^{-1},\omega\rho^{-2}\}\in S(w^{-2}\rho^{-3/2}\la_{\mu}^{-1},\bg),\\
\{\hte,w^{-1}\}\in S(w^{-1}\la_{\mu}^{-1},\bg),\;\{\hte,\omega\rho^{-2}\}\in S(\rho^{-1}\la_{\mu}^{-1},\bg),\;\{w^{-1},\la_{\mu}^{-1}\}\in S(\la_{\mu}^{-1},\bg).
\end{align*}
\end{lem}
 From Lemma \ref{lem:LaM}
it follows that
\begin{equation}
\label{eq:PETIT}
\begin{split}
-2{\mathsf{Im}}(\tLa v,v)
\geq \frac{d}{dx_0}\|v\|^2+\frac{3}{2}\ga\|\la_{\mu}^{\varep}v\|^2
+c\mu n(\chi^2\zeta_{+}^2\hte w^{-1} v,v)\\
+c\mu n\|\chi\zeta_{+}\hte^{1/2}w^{-1/2}v\|^2
-C\|\la_{\mu}^{+0}v\|^2
\end{split}
\end{equation}
with some $c>0$. Let $\zeta_0(s)$, $\chi_0(s)\in C^{\infty}(\R)$ be such that ${\rm supp}\,\zeta_0$ is contained  in $\{\zeta_{+}=1\}$ and $\chi_0=1$ for $s\leq c$ with some $c>0$ and ${\rm supp}\,\chi_0\subset\{\chi=1\}$. Set $\zeta_0=\zeta_0(\hte w^{-1})$ and $\chi_0=\chi_0(\hphi_1^2w^{-1})$. Replace $u$
by $w^{-1}\eta\hte^{1/2} u$, $\eta=\chi_0\zeta_0 $ in  \eqref{eq:PETIT} it follows that
\begin{equation}
\label{eq:kikaku:1}
\begin{split}
-2{\mathsf{Im}}\,(\tLa(w^{-1}\eta \hte^{1/2}u),w^{-1}\eta\hte^{1/2} u)
\geq \frac{d}{dx_0}\|w^{-1}\eta \hte^{1/2} u\|^2\\
+\frac{\ga}{2}\|\la^{\varep}w^{-1}\eta \hte^{1/2} u\|^2
+c\mu n(w^{-1}\chi^2\zeta_{+}^2\hte (w^{-1}\eta\hte^{1/2} u),w^{-1}\eta \hte^{1/2} u).
\end{split}
\end{equation}
We first examine $[\tLa,w^{-1}\eta\hte^{1/2}]u$. Note
$\{\xi_0-\phi_1,w^{-1}\}=-2\hphi_2w^{-3}\{\xi_0-\phi_1,\hphi_2\}$ modulo $S(w^{-1},\bg)$. From \eqref{eq:global:assump} we see $
\eta \hte^{1/2}\{\xi_0-\phi_1,w^{-1}\}-c\eta w^{-2}\hte^{3/2}\in S(w^{-1},\bg)$ with some $c\in S(1,\bg)$. Since $\eta \hte^{1/2}\{\psi,w^{-1}\}=\eta \hte^{1/2}\{{\tilde \zeta}\hte\phi_1,w^{-1}\}$  for  $\chi\chi_2=0$ then noting  Lemma \ref{lem:DGB} we have $
\eta \hte^{1/2}\{\xi_0-\phi_1+\psi,w^{-1}\}-b\eta w^{-2}\hte^{3/2}\in S(w^{-1},\bg)$. 
We next examine $w^{-1}\{\xi_0-\phi_1+\psi,\eta\hte^{1/2}\}$. Since $\hte^{-1/2}\zeta_0\in S(w^{-1/2},\bg)$ from \eqref{eq:global:assump}  we have $w^{-1}\{\xi_0-\phi_1,\eta\hte^{1/2}\}-c\eta\hte w^{-3/2}\in \mu S(w^{-1},\bg)$ by similar arguments. Noting $\{{\tilde\zeta}\hte\phi_1,\eta\hte^{1/2}\}\in \mu S(w^{-1},\bg)$ we get
\begin{equation}
\label{eq:TNT}
\{\tLa,w^{-1}\eta \hte^{1/2}\}-c\eta w^{-2}\hte^{3/2}-c'\eta\hte w^{-3/2}\in  \mu S(w^{-1},\bg)
\end{equation}
with some $c$, $c'\in \mu S(1,\bg)$. From \eqref{eq:TNT} one has
\begin{eqnarray*}
|{\mathsf{Im}}([\tLa,w^{-1}\eta \hte^{1/2}]u,w^{-1}\eta\hte^{1/2} u)|\leq {\mathsf{Re}}(c\eta w^{-2}\hte^{3/2} u,w^{-1}\eta\hte^{1/2} u)\\
+{\mathsf{Re}}(c'\eta\hte w^{-3/2} u,w^{-1}\eta\hte^{1/2} u)
+C\|w^{-1}u\|^2.
\end{eqnarray*}
Write $
{\mathsf{Re}}(w^{-1}\eta\hte^{1/2} \#c\eta w^{-2}\hte^{3/2})={\mathsf{Re}}(w^{-3/2}\eta\hte\# c\#\eta w^{-3/2}\hte)$ 
modulo $S(w^{-2},\bg)$ and $\tLa (w^{-1}\eta\hte^{1/2} u)=w^{-1}\eta \hte^{1/2} (\tLa u)+[\tLa,w^{-1}\eta\hte^{1/2}]u$ we get
\begin{equation}
\label{eq:kikaku:2}
\begin{split}
{\mathsf{Im}}(\tLa w^{-1}(\eta\hte^{1/2} u),w^{-1}\eta\hte^{1/2} u)\leq {\mathsf{Im}}(w^{-1}\eta\hte^{1/2} (\tLa u),w^{-1}\eta\hte^{1/2}u)\\
+C\mu \|\eta w^{-3/2}\hte u\|^2
+C\|\eta w^{-1}\hte^{1/2}u\|^2+C\|w^{-1}u\|^2.
\end{split}
\end{equation}
We now estimate ${\mathsf{Im}}(w^{-1}\eta\hte^{1/2} (\tLa u),w^{-1}\eta\hte^{1/2} u)$. Thanks to Lemma \ref{lem:toyama} one can write
\[
w^{-1}\eta\hte^{1/2}\#w^{-1}\eta\hte^{1/2}=\eta w^{-1/2}\# \eta w^{-3/2}\hte+b\zeta_0w^{-3/2}\hte+R
\]
with $b\in S(1,\bg)$ where $R\in S(w^{-1},\bg)$ and therefore we have
\begin{equation}
\label{eq:kikaku:3}
\begin{split}
{\mathsf{Im}}(w^{-1}\eta\hte^{1/2}(\tLa u),w^{-1}\eta\hte^{1/2} u)\leq 
{\mathsf{Im}}(w^{-1/2}\eta(\tLa u),\eta w^{-3/2}\hte u)
\\
+C(\|\tLa u\|^2+\|\zeta_0w^{-3/2}\hte u\|^2)+C\|w^{-1} u\|^2+C\|u\|^2\\
\leq (c\mu n)^{-1}\|\eta w^{-1/2}\tLa u\|^2+(c\mu n/2)\|\zeta_0 w^{-3/2}\hte u\|^2\\
+C(\|\tLa u\|^2+\|w^{-1} u\|^2)
\end{split}
\end{equation}
where $b\in S(1,\bg)$. Combining \eqref{eq:kikaku:1}, \eqref{eq:kikaku:2} and \eqref{eq:kikaku:3} one obtains
\begin{equation}
\label{eq:kikaku:4} 
\begin{split}
\frac{d}{dx_0}\|w^{-1}\eta\hte^{1/2} u\|^2+\frac{\ga}{2}\|\la^{\varep}w^{-1}\eta \hte^{1/2} u\|^2\\
+c\mu n(w^{-1}\chi^2\zeta^2_{+}\hte (w^{-1}\eta \hte^{1/2} u), w^{-1}\eta \hte^{1/2} u)\\
\leq C(\mu n)^{-1}\|\eta w^{-1/2}\tLa u\|^2+c\mu n\|\zeta_0w^{-3/2}\hte u\|^2\\
+C\big(\|\eta w^{-1}\hte^{1/2}u\|^2+\|w^{-1}u\|^2+\|\tLa u\|^2\big).
\end{split}
\end{equation}
We now estimate $(w^{-1}\chi^2\zeta_{+}^2\hte (w^{-1}\eta \hte^{1/2} u),w^{-1}\eta \hte^{1/2} u)$ from below. Note that 
\[
w^{-1}\eta \hte^{1/2}\#w^{-1}\chi^2\zeta_{+}^2\hte\# w^{-1}\eta \hte^{1/2}=w^{-3}\eta^2\chi^2\zeta^2_{+}\hte^2+R
\]
with $R\in S(w^{-2},\bg)$ and hence we have
\[
(w^{-1}\chi^2\zeta_{+}^2\hte (w^{-1}\eta \hte^{1/2} u),w^{-1}\eta \hte^{1/2} u)\geq (w^{-3}\eta^2\hte^2 u,u)
-C\|w^{-1}u\|^2
\]
for $\eta^2\chi^2\zeta^2_{+}=\eta^2$. Since $\eta^2w^{-3}\hte^2=\eta w^{-3/2}\hte\#\eta w^{-3/2}\hte+R$ with $R \in S(w^{-2},\bg)$ which proves that
\begin{equation}
\label{eq:kikaku:b}
\begin{split}
2{\mathsf{Re}}(w^{-1}\chi^2\zeta_{+}^2\hte(\eta w^{-1}\hte^{1/2}u),\eta w^{-1}\hte^{1/2} u)\geq (\eta^2w^{-3}\hte^2 u,u)\\
+\|\eta w^{-3/2}\hte u\|^2-C\|w^{-1}u\|^2.
\end{split}
\end{equation}
 Write $
\zeta_0^2w^{-3}\hte^2=\eta^2w^{-3}\hte^2+(1-\chi_0^2)\zeta_0^2w^{-3}\hte^2$ and consider 
\[
\mu^{-2}M\zeta_0^2\hte\phi_1^2-(1-\chi_0^2)\zeta_0^2w^{-3}\hte^2
=(\mu^{-1}H\zeta_0\hte^{1/2}\phi_1)^2
\]
where $H=(M-(1-\chi_0^2)w^{-3}\la_{\mu}^{-2}\hphi_1^{-2}\hte)^{1/2}\in S(1,\bg)$ for large $M>0$ which follows from $(1-\chi_0^2)\hphi_1^{-2}\in S(w^{-1},\bg)$. Then it is not difficult to see
\begin{align*}
 (\mu^{-1}H\zeta_0\hte^{1/2}\phi_1)\#(\mu^{-1}H\zeta_0\hte^{1/2}\phi_1)=(\mu^{-1}H\zeta_0\hte^{1/2}\phi_1)^2\\
 +b_1w^{-1/2}\phi_1+b_2w^{-1}\hphi_1^2\la+R
 \end{align*}
with $b\in \mu^{-1}S(1,\bg)$ and $R\in S(w^{-2},\bg)$.
 Noting $b_1w^{-1/2}\phi_1=w^{-1}\#b_1w^{1/2}\phi_1+R_1$ and $b_2w^{-1}\hphi_1^2\la=w^{-1}\#b_2\hphi_1^2\la+R_2$ with $R_i\in S(w^{-2},\bg)$ we conclude that  $
 \mu^{-2}M(\zeta_0^2\hte\phi_1^2 u,u)\geq ((1-\chi_0^2)\zeta_0^2 w^{-3}\hte^2 u,u)$ modulo a term $
C\mu^{-2}(\|b_1w^{1/2}\phi_1 u\|^2 +\|b_2\hphi_1^2\la u\|^2+\|w^{-1}u\|^2)$ which proves together with \eqref{eq:kikaku:b}
\begin{align*}
|(\zeta_0^2w^{-3}\hte^2 u,u)|\leq C\mu^{-2}(a^2{\tilde \zeta}\hte\phi_1^2 u,u)
+2(w^{-1}\chi^2\zeta_{+}^2\hte(\eta w^{-1}\hte^{1/2})u,\eta w^{-1}\hte^{1/2}u)\\
+C\mu^{-2}\big(\|w^{1/2}\phi_1u\|^2+\|\hphi_1^2\la u\|^2+\|w^{-1}u\|^2\big)
\end{align*}
with some $C>0$. To simplify notations we introduce
\begin{defin}
\label{dfn:gosa}\rm  We denote by $O(E)$ a symbol or the set of symbols of the form
\begin{align*}
a_1\mu w^{-1}+a_2\mu \omega^{-1}+a_3\mu^{1/2}\la^{1/2}+a_4w^{1/2}\phi_1\\
+a_5\omega^{1/2}\phi_1
+a_6\phi_2
+a_7\hphi_1^2\la+a_8w\la+a_9\omega\la
\end{align*}
with $a_i\in S(\la_{\mu}^{+0},\bg)$. We denote $S(\la^{t_1}\la_{\mu}^{t_2}w^s,\bg)O(E)$ a symbol or the set of symbols which is a linear combination of $w^{-1}$, $\omega^{-1}$, $\la^{1/2}$, $w^{1/2}\phi_1$, $\omega^{1/2}\phi_1$, $\phi_2$, $\la\hphi_1^2$, $w\la$ and $\omega\la$ with coefficients in $S(\la^{t_1}\la_{\mu}^{t_2}w^s,\bg)$. We also denote 
\begin{align*}
\|O(E)u\|^2=\mu^2(\|w^{-1}u\|^2+\|\omega^{-1}u\|^2)+\mu\|\la^{1/2}u\|^2+\|w^{1/2}\phi_1u\|^2\\
+\|\omega^{1/2}\phi_1u\|^2+\|\phi_2u\|^2+\|\hphi_1^2\la u\|^2+\|w\la u\|^2+\|\omega\la u\|^2.
\end{align*}
\end{defin}
\begin{pro}
\label{pro:wwLa} Let $\chi_0$, $\zeta_0$ be as above. Then we have
\begin{align*}
C\mu n\|\chi_0\zeta_0w^{-1/2}\tLa u\|^2+C\mu^2n^2\|\tLa u\|^2
\geq c\mu^2n^2\frac{d}{dx_0}\|\chi_0\zeta_0 w^{-1}\hte^{1/2} u\|^2\\
+c\mu^2n^2\|\chi_0\zeta_0 w^{-1}\hte^{1/2}\la^{\varep} u\|^2
+c\mu^3n^3(\|\zeta_0w^{-3/2}\hte u\|^2
+(\zeta_0^2w^{-3}\hte^2 u,u))\\
-C\mu({\tilde \zeta}a^2\hte \phi_1^2 u,u)-C\mu \|O(E)u\|^2
\end{align*}
with some $c>0$ and $C=C(n)$.
\end{pro}
Replacing $u$ now 
by $w^{-1/2}\eta\hte^{1/2} u$ in  \eqref{eq:PETIT}
 and repeating similar arguments we obtain
\begin{pro}
\label{pro:wwLam} Let $\chi_0$, $\zeta_0$ be as above. Then we have
\begin{align*}
C\|\tLa u\|^2\geq c\mu n\frac{d}{dx_0}\|\chi_0\zeta_0 w^{-1/2}\hte^{1/2} u\|^2
+c\mu n\ga(\|\chi_0\zeta_0 w^{-1/2}\hte^{1/2}\la^{\varep} u\|^2\\
+c\mu^2n^2(\|\zeta_0w^{-1}\hte u\|^2
+(\zeta_0^2w^{-2}\hte^2 u,u))
-C\gamma^{1/2}\|O(E)u\|^2.
\end{align*}
\end{pro}
%


\section{Transformed symbol ${\tilde Q}$}

We start with
\begin{lem}
\label{lem:siroebi}
One can write $
O(E)=T\#(O(E)+R)$ 
with $R\in S(\la_{\mu}^{+0},\bg)$.
\end{lem}
%
\noindent
Proof:
 Let $A\in O(E)$. Then it is easy to check that $T^{(\al)}_{(\be)}A^{(\be)}_{(\al)}\in S(\la_{\mu}^{-1/4},\bg)O(E)$ for $|\al+\be|=1$. Then we have $TA-T\#A=TA_1$ with $A_1\in S(\la_{\mu}^{-1/4},\bg)O(E)$. Repeating the same arguments we get $TA=T\#(A+A_1+\cdots +A_4)+K$ where $K\in S(\la_{\mu}^{-1},\bg)O(E)\subset S(\la_{\mu}^{+0},\bg)$. Since $T\#T^{-1}=1-r$ with $r\in \mu^{1/4}S(1,\bg)$ and hence the inverse of $1-r$ exists in ${\mathcal L}(L^2,L^2)$ which is given by $Op(b)$ with $b\in S(1,\bg)$ (see \cite{Bea}) and hence $T\#{\tilde T}=1$ with ${\tilde T}=T^{-1}\#b\in S(\la_{\mu}^{+0},\bg)$. Then writing $K=T\#({\tilde T}\#K)$ we get the assertion.
\qed
 
Recall $W^{\al}_{\be}=T^{-1}\dif_x^{\be}\dif_{\xi'}^{\al}T\in S(\la_{\mu}^{-3|\al|/4+|\be|/2+0},\bg)$. Since $q^{(\al)}_{(\be)}\in S(\la^2\la_{\mu}^{-|\al|},\bg)$ for $|\al+\be|=1$ by Lemma \ref{lem:toyama} we see 
\begin{equation}
\label{eq:nodo}
\begin{split}
q\#T=T\#q-i nT\{q,\Psi\}
+\frac{i}{8}T\sum_{|\al+\be|=3}\frac{(-1)^{|\be|}}{\al!\be!}\big(q^{(\al)}_{(\be)}W^{\be}_{\al}-W^{\al}_{\be}q^{(\be)}_{(\al)}\big)
\end{split}
\end{equation}
modulo $R\in \mu^{3/2} S(\la^{1/2+0},\bg)$. We first check
\begin{lem}
\label{lem:q:one} We have 
\[
\sum_{|\al+\be|=3}(-1)^{|\be|}(q^{(\al)}_{(\be)}W^{\be}_{\al}-W^{\al}_{\be}q^{(\be)}_{(\al)})/\al!\be!\in \mu S(\la_{\mu}^{+0},\bg) O(E).
\]
\end{lem}
%
\noindent
Proof:  Write $q=\phi_2^2+f\phi_1^2$ with $f=2a^2{\tilde\zeta}\hte+2a^2\chi_2\hphi_1^2$ and recall $f\in S(1,g)$ and $f^{(\al)}_{(\be)}\in S(\la_{\mu}^{-|\al|},g)$ for $|\al+\be|=1$ by Lemma \ref{lem:toyama}. Applying Lemma \ref{lem:toyama} again one can check that $(f\phi_1^2)^{(\al)}_{(\be)}W^{\be}_{\al}$ with $|\al+\be|=3$ is a linear combination of $\la^{1/2}$, $w^{1/2}\phi_1$ and $\hphi_1^2\la$ with coefficients in $\mu S(\la_{\mu}^{+0},\bg)$ which proves the assertion. 
\qed
 
 We make more detailed studies on $\{q,\Psi\}$ and $\{\Psi,\{q,\Psi\}\}$. Let us denote $\Psi_1=\zeta^2\chi^2\log{\phi}\in S(\la_{\mu}^{+0},g)$ and $\Psi_2=\zeta^2\Phi$, $\Phi\in S(1,g_1)$ so that $\Psi=\Psi_1+\Psi_2$.   
\begin{lem}
\label{lem:simizu}
We have
\begin{align*}
\{q,\Psi\}=\nu\zeta_{+}^2a^2\hte\phi_1(\Gamma+\Delta)\{\phi_1,\hphi_2\}
+a_1 O(E)+a_2O(E)+a_3O(E)
\end{align*}
where $a_1=\zeta\chi a_1'$ with $a'_1\in \mu S(w^{-1/2}\la_{\mu}^{+0},\bg)$ and $a_2=\zeta a_2'$ with $a'_2\in \mu S(\rho^{-1/2},\bg)$  and $a_3\in \mu S(\la_{\mu}^{+0},\bg)$.
\end{lem}
%
\noindent
Proof:
Thanks to \eqref{eq:Fai} and Lemma \ref{lem:DGB} we can see $
\{\phi_2^2+a^2\chi_2\hphi_1^4\la^2,\Psi_2\}=a_2O(E)$ 
where $a_2=\zeta a_2'$ with $a_2\in\mu S(\rho^{-1/2},\bg) $. Similarly from Lemma \ref{lem:DGB} and 
\begin{equation}
\label{eq:soji}
\{F,\log{\phi}\}=\{F,\hphi_2\}/w+\{F,\la_{\mu}^{-1}\}/(2w\phi)
\end{equation}
we obtain $\{\phi_2^2+a^2\chi_2\hphi_1^4\la^2,\Psi_1\}=a_1O(E)$ with $a_1\in \mu S(w^{-1/2}\la_{\mu}^{+0},\bg)$ where clearly $a_1=\zeta\chi a_1'$. We turn to $\{a^2{\tilde\zeta}\hte\phi_1^2,\Psi_j\}$. Repeating similar arguments one can check that $\{a^2{\tilde\zeta}\hte\phi_1^2,\Psi_j\}=a^2{\tilde\zeta}\hte\{\phi_1^2,\Psi_j\}+a_jO(E)$ where $a_j$ verifies the same properties as above. From the same arguments proving \eqref{eq:nisiyama} and \eqref{eq:arata:b} one can show
\[
\{\phi_1,\Psi\}=(\zeta^2\Gamma+\Delta)\{\phi_1,\hphi_2\}+a_3O(E)
\]
with $a_3\in \mu S(\la_{\mu}^{+0},\bg)$. Since ${\tilde\zeta}\zeta=\nu\zeta_{+}^2$ and ${\tilde\zeta}\Delta=\nu\zeta_{+}^2\Delta$ we get the assertion.
\qed
\begin{lem}
\label{lem:tanaka}
We have
\begin{align*}
\{\Psi,\{q,\Psi\}\}=\zeta_{+}(a_{1}\mu w^{-2}\hte+a_2w^{-1/2}\hte\phi_1+a_3\phi_1)+\mu S(\la_{\mu}^{+0},\bg)O(E)
\end{align*}
where $a_{j}\in \mu S(\la_{\mu}^{+0},\bg)$.
\end{lem}
%
\noindent
Proof:
By obvious abbreviated notation  we see $\Psi^{(\al)}_{(\be)}(a_jO(E))^{(\be)}_{(\al)}\in O(E)$ for $|\al+\be|=1$ and hence $\{\Psi,\sum a_jO(E)\}\in O(E)$. 
With $b=\zeta_{+}^2a^2\hte\Delta\{\phi_1,\hphi_2\}\in \mu S(\la_{\mu}^{+0},\bg)$ noting $\hte\Delta\in S(\la_{\mu}^{+0},\bg)$ it is easy to check that $\{\Psi, b\phi_1\}$ is a linear combination of $w^{1/2}\phi_1$ and $w\la$ with coefficients $\mu S(\la_{\mu}^{+0},\bg)$ because $\zeta^2\Phi\in S(1,{\tilde g})$. Therefore $
\{\Psi,b\phi_1\}\in \mu S(\la_{\mu}^{+0},\bg)O(E)$. 
Recall $\Gamma=r+2\omega\rho^{-2}$ and note ${\rm supp}\,r\subset {\rm supp}\,\chi$.  With $B_1=\zeta_{+}^2a^2\hte\phi_1 r\{\phi_1,\hphi_2\}$ taking Lemma \ref{lem:DGB} into account we can prove that
\[
\{\Psi,B_1\}=\zeta_{+}(a_{1}w^{-1/2}\hte\phi_1+a_{2}\phi_1+a_{3}\mu w^{-2}\hte)
\]
where $a_{j}\in \mu S(\la_{\mu}^{+0},\bg)$. Here it is obvious that the supports of $a_{j}$ are contained in that of $\zeta_{+}$. By similar arguments we get with $B_2=2\zeta_{+}^2a^2\hte\phi_1 \omega\rho^{-2}\{\phi_1,\hphi_2\}$
\[
\{\Psi,B_2\}=\zeta_{+}({\tilde a}_{1}w^{-1/2}\hte\phi_1+{\tilde a}_{2}\phi_1+{\tilde a}_{3}\mu w^{-2}\hte)
\]
where ${\tilde a}_{j}\in \mu S(\la_{\mu}^{+0},\bg)$. This proves the assertion.
\qed
\begin{pro}
\label{pro:jiten} We have
\begin{align*}
q\#T=T\#\big(q-in\{q,\Psi\}-n^2\{\Psi,\{\Psi,q\}\}+i(
a_1\mu w^{-3/2}\hte\\
+a_2w^{1/2}\hte\la+a_3\hte\phi_1)+\mu S(\la_{\mu}^{+0},\bg)O(E)
\end{align*}
where $a_j\in \mu S(\la_{\mu}^{+0},\bg)$ are real valued and ${\rm supp}\,a_j\subset {\rm supp}\,\zeta_{+}$. 
\end{pro}
%
\noindent
Proof:
 From Lemma \ref{lem:simizu}  it is clear that $T^{-1}T^{(\al)}_{(\be)}\{q,\Psi\}^{(\be)}_{(\al)}$ is $c_1\hte\phi_1+c_2w^{1/2}\hte\la$ modulo $\mu^{1/2} S(\la_{\mu}^{+0},\bg)O(E)$ because $\la_{\mu}^{-1/4}\in S(w^{1/2},\bg)$ for $|\al+\be|=2$. Therefore we get
\begin{align*}
T\#\{q,\Psi\}=T\{q,\Psi\}+nT\{\Psi,\{q,\Psi\}\}/2i+c_1\hte\phi_1\\
+c_2 w^{1/2}\hte\la 
+\mu S(\la_{\mu}^{+0},\bg)O(E)
\end{align*}
where $c_i\in \mu S(\la_{\mu}^{+0},\bg)$ is real. It is clear that ${\rm supp}\,c_j\subset {\rm supp}\,\zeta_{+}$. Thus we have $
T\{q,\Psi\}=T\#\big(\{q,\Psi\}+n\{\Psi,\{\Psi,q\}\}/2i-c_1\hte\phi_1 -c_2\mu w^{1/2}\hte\la\big) +\mu S(\la_{\mu}^{+0},\bg)O(E)$.
From Lemma \ref{lem:tanaka} it can be seen that $\Psi^{(\al)}_{(\be)}\{\Psi,\{q,\Psi\}\}^{(\be)}_{(\al)}$ for $|\al+\be|=1$ are written as $
a_1\mu w^{-3/2}\hte+a_2\omega^{1/2}\hte\la+c_2\hte\phi_1$ modulo $\mu S(\la_{\mu}^{+0},\bg)O(E)$. This proves 
\begin{align*}
T\{q,\Psi\}=T\#(\{q,\Psi\}+n\{\Psi,\{\Psi,q\}\}/2i+a_1\mu w^{-3/2}\hte\\
+a_2w^{1/2}\hte\la+a_3\hte\phi_1)+\mu S(\la_{\mu}^{+0},\bg)O(E)
\end{align*}
where $a_j\in \mu S(\la_{\mu}^{+0},\bg)$ with ${\rm supp}\,a_j\subset {\rm supp}\,\zeta_{+}$. This proves the assertion. 
\qed
 
Taking $\hte \Delta\in S(\la_{\mu}^{+0},\bg)$ into account we have

\begin{cor}
\label{cor:QIm} We have
\begin{align*}
{\mathsf{Im}}\,{\tilde Q}=T_2-\nu n\zeta_{+}^2a^2\hte\phi_1 \Gamma \{\phi_1,\hphi_2\}+a_1\mu w^{-3/2}\hte
+a_2w^{1/2}\hte\la \\
+a_3\hte\phi_1
+c_1O(E)+c_2O(E)+\mu S(\la_{\mu}^{+0},\bg)O(E)
\end{align*}
where $a_j \in \mu S(\la_{\mu}^{+0},\bg)$ are real valued of which support is contained in ${\rm supp}\,\zeta_{+}$ and $c_1=\zeta\chi c'_1$ with $c'_1\in \mu S(w^{-1/2}\la_{\mu}^{+0},\bg)$  and $c_2=\zeta c_2'$ with $c_2\in \mu S(\rho^{-1/2},\bg)$.
\end{cor}
\begin{cor}
\label{cor:QRe} We have
\begin{align*}
{\mathsf{Re}}\,{\tilde Q}=q+T_1+\zeta_{+}(a_{1}\mu w^{-2}\hte+a_2w^{-1/2}\hte\phi_1
+a_3\phi_1)+\mu S(\la_{\mu}^{+0},\bg)O(E)
\end{align*}
where $a_{j}\in \mu S(\la_{\mu}^{+0},\bg)$.
\end{cor}
\section{Estimate $(({\mathsf{Re}}\,{\tilde Q}-T_1+{\bar\kappa}\mu\la)u,u)$}

Here we write $q+T_1=q+{\bar\kappa}\mu\la+(T_1-{\bar\kappa}\mu\la)$ and instead of $q$ we consider $q+{\bar\kappa}\mu\la$ with ${\bar\kappa}>0$ in \eqref{eq:mao}. 
 In this section we study $(({\mathsf{Re}}\,{\tilde Q}-T_1+{\bar\kappa}\mu \la)u,u)$. Without restrictions we can assume ${\bar\kappa} =1$.
\begin{pro}
\label{pro:qandQ} Let   $c_{\pm}\in S(1,\bg)$ be real. Then 
we have
\begin{equation}
\label{eq:Qnojitu}
\begin{split}
C\big((q+\mu\la)u,u)\geq \sum (\|c_{\pm}\zeta_{\pm}|\hte|^{1/2}\phi_1 u\|^2
+|(c_{\pm}\zeta_{\pm}^2|\hte|\phi_1^2 u,u)|)\\
+|(\phi_2^2u,u)|+|(w\phi_1^2 u,u)|+|(\omega\phi_1^2u,u)|\\
+|(w^2\la^2 u,u)|
+|(\omega^2\la^2 u,u)|+\|O(E)u\|^2.
\end{split}
\end{equation}
\end{pro}
%
\noindent
Proof:
One can write
\begin{align*}
Ma^2{\tilde \zeta}\hte\phi_1^2-(c_{+} \zeta_{+}^2\hte +c_{-}\zeta_{-}^2|\hte|)\phi_1^2
=H_{+}^2\zeta_{+}^2\hte\phi_1^2+H_{-}^2\zeta_{-}^2|\hte|\phi_1^2
\end{align*}
with $H_{+}=(Ma^2\nu-c_{+})^{1/2}$ and $H_{-}=(Ma^2{\hat h}-c_{-})^{1/2}$ where $M>0$ is chosen so that $Ma^2\nu-c_{+}\geq c$, $Ma^2{\hat h}-c_{-}\geq c>0$. Since $\zeta_{\pm}|\hte|^{1/2}\in S(|\hte|^{1/2}, \bg)$ by Lemma \ref{lem:ameagari} noting  $H_{\pm}\in S(1,\bg)$ we can write
\begin{align*}
\zeta_{\pm}|\hte|^{1/2}\phi_1H_{\pm}\#\zeta_{\pm}|\hte|^{1/2}\phi_1H_{\pm}-\zeta_{\pm}^2|\hte|\phi_1^2H_{\pm}^2\\
=\sum_{|\al+\be|=2}C_{\al\be}(\zeta_{\pm}|\hte|^{1/2}\phi_1H_{\pm})^{(\al)}_{(\be)}(\zeta_{\pm}|\hte|^{1/2}\phi_1H_{\pm})_{(\al)}^{(\be)}\\=b_1w^{-3}\hphi_1^2+b_2w^{-5/2}\hphi_1
\end{align*}
modulo $\mu^2 S(w^{-2},\bg)$ where $b_i\in \mu^2 S(1,\bg)$. Write $
b_1w^{-3}\hphi_1^2=c_1w^{1/2}\phi_1\#w^{1/2}\phi_1+R_1$ and $b_2w^{-5/2}\hphi_1=\mu c_2w^{-1}\#w^{1/2}\phi_1+R_2$ with $c_i\in S(1,\bg)$ and $R_i\in \mu S(w^{-2},\bg)$ 
we conclude 
\begin{equation}
\label{eq:kitune}
\sum|(c_{\pm}\zeta_{\pm}^2|\hte|\phi_1^2u,u)|\leq M(a^2{\tilde\zeta}\hte\phi_1^2 u,u)+C(\|w^{1/2}\phi_1u\|^2+\mu^2\|w^{-1}u\|^2).
\end{equation}
Similarly $c_{\pm}^2\zeta_{\pm}^2|\hte|\phi_1^2$ can be written
\[
c_{\pm}\zeta_{\pm}|\hte|^{1/2}\phi_1\#c_{\pm}\zeta_{\pm}|\hte|^{1/2}\phi_1+b_{\pm}\mu^2 w^{-5/2}\hphi_1+b'_{\pm}\mu^2w^{-3}\hphi_1^2+R
\]
with $R\in S(w^{-2},\bg)$. Thus $\|c_{\pm}\zeta_{\pm}|\hte|^{1/2}\phi_1 u\|^2$ is estimated also by the right-hand side of \eqref{eq:kitune}.

We next study ${\tilde q}=\phi_2^2+\chi_2a^2\phi_1^4\la^{-2}+\mu\la$. If $\chi_2\neq 1$ so that $\hphi_1^2\leq d_3w$ it is clear $\hphi_1^4\leq C(\hphi_2^2+\la_{\mu}^{-1})$ and then noting $\la\la_{\mu}^{-1}=\mu$ we have  ${\tilde q}\geq c\hphi_1^4\la^2$ with some $c>0$. If $\chi_2=1$ this inequality is obvious. Since $\phi_1^4\la^{-2}+\mu\la=\la^2\omega^2$ and $\phi_2^2+\mu\la=w^2\la^2$ it is obvious ${\tilde q}\geq c(w^2+\omega^2)\la^2$ with some $c>0$. Let us set ${\tilde q}-c\omega^2\la^2=F^2$ with $F=\la({\tilde q}\la^{-2}-c\omega^2)^{1/2}\in S(\la,\bg)$. If we note $\chi_2a^2\hphi_1^4\in S(w^2,g)$ and $\omega\in S(\omega, G_1)$ with $G_1=\omega^{-1/2}|dx|^2+\omega^{-1/2}\xim^{-2}|d\xi'|^2$ then it is not difficult to see that $
F^2=F\#F+R$ with $R\in \mu^2 S(w^{-1}+\omega^{-1},\bg)$. Thus  we conclude that
\begin{align*}
({\tilde q} u,u)\geq c(\omega^2\la^2 u,u)
-C\mu^4(\|\omega^{-1}u\|^2+\|w^{-1}u\|^2)-C\|u\|^2.
\end{align*}
Repeating a similar argument we get $({\tilde q}u,u)\geq c\|\hphi_1^2\la u\|^2-C\mu^4\|w^{-1}u\|^2-C\|u\|^2$.
Since $\omega^2\la^2=\omega\la\#\omega\la+R$ with $R\in \mu^2 S(\omega^{-1},\bg)$ and hence $(\omega^2\la^2 u,u)\geq \|\omega\la u\|^2-C(\mu^4\|\omega^{-1}u\|^2+\|u\|^2)$. 
 Recalling  $\phi_2^2+\mu\la=w^2\la^2$ similar arguments show
\[
((\phi_2^2+\mu\la)u,u)\geq c (w^2\la^2u,u)+\|w\la u\|^2-C(\mu^4\|w^{-1}u\|^2+\|u\|^2).
\]
Noting $\mu\la\leq  w^2\la^2 \in S(\la^2,g_0)$ we see $(w^2\la^2 u,u)\geq \mu\|\la^{1/2}u\|^2-C\|u\|^2$. On the other hand since one can write $w^{-1}=(w^{-1}\la^{-1/2})\#\la^{1/2}+R$ with $R\in S(1,\bg)$ remarking $w^{-1}\la^{-1/2}\in \mu^{-1/2}S(1,\bg)$ we have $\|w^{-1}u\|^2\leq C\mu^{-1}\|\la^{1/2}u\|^2+C\|u\|^2$. Similarly we have $\|\omega^{-1}u\|^2\leq C\mu^{-1}\|\la^{1/2}u\|^2+C\|u\|^2$. Thus we get 
\begin{equation}
\label{eq:Wtwo}
\begin{split}
\mu^2(\|w^{-1}u\|^2+\|\omega^{-1}u\|^2)+\mu\|\la^{1/2}u\|^2+\|w\la u\|^2+\|\omega\la u\|^2\\
+|(w^2\la^2 u,u)|+|(\omega^2\la u,u)|+\|\hphi_1^2\la u\|^2
\leq C({\tilde q}u,u)+C\|u\|^2.
\end{split}
\end{equation}
Note $
w^{1/2}\phi_1\#w^{1/2}\phi_1=w\phi_1^2+R$ with $R\in \mu^2 S(w^{-2},\bg)$ and $
w\phi_1^2={\mathsf{Re}}(\la\hphi_1^2\# w\la )+R$ with $R\in \mu^2 S(w^{-2},\bg)$ 
we have
\[
\|w^{1/2}\phi_1 u\|^2+|(w\phi_1^2 u,u)|\leq C\|O(E)u\|^2.
\]
We get $\|\omega^{1/2}\phi_1 u\|^2+|(\omega\phi_1^2 u,u)|\leq C(\|\la\hphi_1^2 u\|^2+\|\omega\la u\|^2+\mu^2\|\omega^{-1}u\|^2)$ by a repetition of similar arguments.
It is easy to see $
\|\phi_2 u\|^2+|(\phi_2^2u,u)|\leq C((\phi_2^2+\mu\la)u,u)+\|u\|^2)$ then we conclude the assertion by \eqref{eq:Wtwo}.
\qed
%
\begin{cor}
\label{cor:kitagata} We have $
\|\hte \phi_1 u\|^2+|( \hte \phi^2_1 u,u)|\leq C((q+\mu\la)u,u)+C\|u\|^2$. 
\end{cor}
%
\noindent
Proof:
 Take $\eta(s)\in C_0^{\infty}(\R)$ so that $\zeta_{-}+\zeta_{+}+\eta=1$. Thanks to Proposition \ref{pro:qandQ} it suffices to prove $|(\eta \hte \phi^2_1 u,u)|\leq C((q+\mu\la)u,u)+C\|u\|^2$. 
Note that one can write $
\eta\hte\phi_1^2=cw\phi_1^2$ 
then the assertion follows immediately.
%
\qed
\begin{lem} 
\label{lem:tango} Let $\chi_0=\chi_0(\hphi_1^2w^{-1})$ with $\chi_0(s)\in C_0^{\infty}(\R)$ which is $1$ near $s=0$. Then we have 
\[
((1-\chi_0)\zeta_{\pm}^2|\hte| w\la^2 u,u)\leq C((q+\mu\la)u,u)+C\|u\|^2. 
\]
\end{lem}
%
\noindent
Proof:
Note that $
Ma^2{\tilde \zeta}\hte\phi_1^2-(1-\chi_0)( \zeta_{+}^2\hte +\zeta_{-}^2|\hte|)w\la^2
=H_{+}^2\zeta_{+}^2\hte\phi_1^2+H_{-}^2\zeta_{-}^2|\hte|\phi_1^2$ 
where $H_{+}=(Ma^2\nu-(1-\chi_0)w\hphi_1^{-2})^{1/2}$ and $H_{-}=(Ma^2{\hat h}-(1-\chi_0)w\hphi_1^{-2})^{1/2}$
 which are in $S(1,\bg)$ taking $M>0$ large. The rest of the proof is just a repetition of the proof of Proposition \ref{pro:qandQ}.
\qed
 
It is easy to check 
\begin{equation}
\label{eq:seki}
\begin{split}|(\zeta_{+}(a_{1}\mu w^{-2}\hte+a_2w^{-1/2}\hte\phi_1
+a_3\phi_1)u, u)|\\
\leq Cn(\mu^2+\ga^{-1/2})\|\zeta_{+}w^{-1}\hte u\|^2
+Cn^{-1}\|O(E)u\|^2.
\end{split}
\end{equation}
From Propositions \ref{pro:qandQ} and \ref{pro:wwLam}  and Corollary \ref{cor:kitagata} together with \eqref{eq:seki}  we obtain
\begin{pro}
\label{pro:kasitu} There exist $\ga_0>0$, $\mu_0>0$, $n_0>0$ such that we have
\[
C(({\mathsf{Re}}\,{\tilde Q}-T_1+{\bar\kappa}\mu \la)u,u)+C\|\tLa u\|^2\geq |(\hte\phi^2_1 u,u)|+\|\hte\phi_1 u\|^2+\|O(E)u\|^2.
\]
for $\ga\geq \ga_0$, $0<\mu<\mu_0$ and $n\geq n_0$. We have also
\begin{align*}
C(\la_{\mu}^{2\varepsilon}({\mathsf{Re}}\,{\tilde Q}-T_1+{\bar\kappa}\mu \la)u,u)+C\|\la_{\mu}^{\varepsilon}\tLa u\|^2
\geq\|\la_{\mu}^{\varepsilon}\hte\phi_1 u\|^2+\|\la_{\mu}^{\varepsilon}O(E)u\|^2.
\end{align*}
\end{pro}
%

\section{ Estimate ${\mathsf{Re}}(({\mathsf{Re}}\, {\tilde Q}-T_1+{\bar\kappa}\mu\la)u,({\mathsf{Im}}\,{\tilde \la})u)$}

Recall Lemma \ref{lem:whatla} which gives $
{\mathsf{Im}}{\tilde \la}=n {\tilde e}_1\Gamma\zeta^2_{+}\hte+R_1$ with $R_1\in S(\la_{\mu}^{+0},\bg)$.
Denote ${\tilde q}=\phi_2^2+\chi_2a^2\phi_1^4\la^{-2}+\mu\la$ again. Note ${\mathsf{Re}}({\tilde e}_1\hte\zeta_{+}^2 \Gamma\#{\tilde q})={\tilde e}_1\zeta_{+}^2\hte {\tilde q} \Gamma+R$ 
with $R\in \mu S(\la^{1+0},\bg)$ since $ \Gamma\in S(w^{-1}\la_{\mu}^{+0},{\tilde g})$ and $\phi_2^2+\chi_2a^2\phi_1^4\la^{-2}\in S(w^2\la^2,{\tilde g})$. Thus noting $|(Ru,u)|\leq C\mu\|\la^{1/2+0}u\|^2)$ we get
\begin{align*}
{\mathsf{Re}}({\tilde q}u,{\tilde e}_1\hte\zeta_{+}^2\Gamma u)\geq ({\tilde e}_1\zeta_{+}^2\hte {\tilde q}\Gamma u,u)
-C\|\la_{\mu}^{+0}O(E)u\|^2.
\end{align*}
Write $
M{\tilde e}_1\zeta_{+}^2\hte {\tilde q} \Gamma-\mu\zeta_{+}^2\hte w^2\la^2 \Gamma=H\#(M{\tilde e}_1{\tilde q}w^{-2}\la^{-2}-\mu)\Gamma\#H+R$ 
with $H=\zeta_{+}\hte^{1/2}w\la$ and $R\in \mu  S(\la ^{1+0},\bg)$. Since $0\leq (M{\tilde e}_1{\tilde q}w^{-2}\la^{-2}-\mu) \Gamma\in \mu S(w^{-1}\la_{\mu}^{+0},\bg)\subset \mu S^{1/2+0}_{3/4,1/2}$ then from the Fefferman-Phong inequality it follows that
\begin{align*}
M({\tilde e}_1\zeta_{+}^2\hte {\tilde q} \Gamma u,u)-\mu (\zeta_{+}^2\hte w^2\la^2 \Gamma u,u)
\geq -C\|\la_{\mu}^{+0}O(E)u\|^2.
\end{align*}
Since $
\zeta_{+}^2\hte w^2\la^2 \Gamma=H\#\Gamma\# H+R$ 
with $R\in \mu S(\la^{1+0},\bg)$ taking $(\zeta_{+}^2-\zeta^2)\hte\in S(w,g)$ into account we conclude
\begin{equation}
\label{eq:ReQImla:a} 
\begin{split}
\mu ( \Gamma(\zeta\hte^{1/2}w\la u),\zeta_{+}\hte^{1/2}w\la u)+\mu |(\zeta^2\hte w^2\la^2 \Gamma u,u)|\\
\leq M{\mathsf{Re}}({\tilde q}u,{\tilde e}_1\hte\zeta_{+}^2 \Gamma u)
+C\|\la_{\mu}^{+0}O(E)u\|^2.
\end{split}
\end{equation}
Noticing $\Gamma=r+2\omega\rho^{-2}$ and $\omega^s r\in S(w^{s-1}\la_{\mu}^{+0},{\tilde g})$ for $s\geq 0$ a repetition of similar argument for $\omega$ instead of $w$ shows \eqref{eq:ReQImla:a} where $w$ is replaced by $\omega$. It is easy to check that $w+\omega^3\rho^{-2}\geq c\rho$ with some $c>0$. Since $(w^2+\omega^2)\Gamma\geq \chi^2w+\omega^3\rho^{-2}$ and $C\omega\geq \rho\geq \omega$ on the support of $1-\chi^2$ we see easily that
\[
C\rho\geq (w^2+\omega^2)\Gamma\geq c\rho
\]
with some $c>0$. Then applying the Fefferman-Phong inequality one obtains $(\zeta^2\hte(w^2+\omega^2)\la^2\Gamma u,u)\geq c\|\zeta\hte^{1/2}\rho^{1/2}\la u\|^2-C\|\la_{\mu}^{+0}O(E)u\|^2$. Thus 
\begin{lem}
\label{lem:JCB} We have
\begin{align*}
\mu|(\zeta^2\hte \rho \la^2 u,u)|+\mu\|\zeta\hte^{1/2}\rho^{1/2}\la  u\|^2
\leq C{\mathsf{Re}}({\tilde q}u,{\tilde e}_1\hte\zeta_{+}^2 \Gamma u)+C\|\la_{\mu}^{+0}0(E)u\|^2.
\end{align*}
\end{lem}

We turn to ${\mathsf{Re}}(a^2{\tilde \zeta}\hte\phi_1^2u, {\tilde e}_1\zeta_{+}^2\hte \Gamma u)$. Since $\Gamma=r+2\omega\rho^{-2}$ and $r\hte^2\in S(w,\bg)$ and $\omega\rho^{-2}\hphi_1^2\in S(1,\bg)$ we see that ${\mathsf{Re}}({\tilde e}_1\zeta_{+}^2\hte \Gamma\#a^2{\tilde\zeta}\hte\phi_1^2)$ is 
\begin{eqnarray*}
\nu {\tilde e}_1\zeta_{+}^4\hte^2 a^2\phi_1^2 \Gamma
+\sum_{|\al+\be|=2}\frac{(-1)^{|\be|}}{(2i)^{|\al+\be|}\al!\be!}({\tilde e}_1\zeta_{+}^2\hte \Gamma)^{(\al)}_{(\be)}(a^2{\tilde\zeta}\hte\phi_1^2)^{(\be)}_{(\al)}+R
\end{eqnarray*}
with $R\in \mu S(\la,\bg)$. Consider $
({\tilde e}_1\zeta_{+}^2 \hte)^{(\al_2)}_{(\be_2)}\Gamma^{(\al_1)}_{(\be_1)}(a^2{\tilde\zeta}\hte)^{(\be'')}_{(\al'')}(\phi_1^2)^{(\be')}_{(\al')}
$ 
for $|\al+\be|=2$. By Lemma \ref{lem:toyama} it is not difficult to see that we can write such a term as
\begin{equation}
\label{eq:takuwa}
\begin{split}
\nu c\zeta_{+}^2w\la^2\hte^2+\zeta_{+}(c_{21}w^{-1/2}\la \hte+c_{22}w^{-1}\phi_1\hte
+c_{23}w^{-3/2}\hte\hphi_1^2\la)\\
+(c_{31}w^{-1/2}\phi_1+c_{32}w^{-1} \hphi_1^2\la)
\end{split}
\end{equation}
 with $c\in \mu S(1,\bg)$ and $c_{ij}\in \mu^2 S(1,g)$. One can estimate the last term applying Proposition \ref{pro:qandQ}. The second term can be estimated thanks to Propositions \ref{pro:wwLa} and  \ref{pro:qandQ}. Indeed writing $c_{23}\zeta_{+}\hte\hphi_1^2\la={\mathsf{Re}}(c_{23}\zeta_{+} w^{-3/2}\hte\# \hphi_1^2\la)+R$ with $R\in S(w^{1/2}\la,\bg)$ we have
\begin{align*}
|{\mathsf{Re}}(c_{23}\zeta_{+}w^{-3/2}\hte\hphi_1^2\la u,u)|\leq C\mu^2\ga^{-1/2}\|\zeta_{+}w^{-3/2}\hte
u\|^2+C\ga^{1/2}\|O(E)u\|^2.
\end{align*}
To estimate the first term in \eqref{eq:takuwa} choosing $\nu>0$ small we write $\zeta_{+}^2\hte w\la^2-\nu c\zeta_{+}^2w\hte^2\la^2=H\#H+R$ with $H=\zeta_{+}\hte^{1/2}w^{1/2}\la(1-\nu c\hte)^{1/2}$ and $R\in S(w^2\la^{2+0},\bg)$  and apply Lemma \ref{lem:JCB}. 
We now prove
\begin{lem}
\label{lem:foo} There are $c>0$ and $\nu_0>0$ such that we have 
\begin{equation}
\label{eq:kaze}
\begin{split}
{\mathsf{Re}}(a^2{\tilde \zeta}\hte\phi_1^2u, {\tilde e}_1\zeta_{+}^2\hte \Gamma u)\geq c\nu \mu (\Gamma(\zeta_{+}\hte\phi_1)u,\zeta_{+}\hte\phi_1 u)\\
-C(\mu^3n+\gamma^{-1/2})\|\zeta w^{-3/2}\hte u\|^2
-C(\mu n^{-1}+\gamma^{-1/2})\|\zeta w^{1/2}\hte\la u\|^2\\-C\|\hte\phi_1 u\|^2-C\gamma^{1/2}(\|\zeta w^{-1}\hte u\|^2+\|O(E)u\|^2)
\end{split}
\end{equation}
for $0<\nu\leq \nu_0$.
\end{lem}
%
\noindent
Proof:
It remains to estimate $\nu {\mathsf{Re}}({\tilde e}_1\zeta_{+}^4\hte^2 a^2\phi_1^2 \Gamma u,u)$ from below. Since $(\zeta_{+}^4-\zeta_{+}^2)\hte^2\phi_1^2\Gamma\in S(w\phi_1^2\la_{\mu}^{+0},\bg)$ it suffices to study $\nu {\mathsf{Re}}({\tilde e}_1\zeta_{+}^2\hte^2 a^2\phi_1^2 \Gamma u,u)$.  Note that 
\begin{align*}
{\mathsf{Re}}(\zeta_{+}\hte\phi_1\#{\tilde e}_1a^2\Gamma\#\zeta_{+}\hte\phi_1)={\tilde e}_1\zeta_{+}^2\hte^2 a^2\phi_1^2 \Gamma-\sum\frac{(-1)^{|\be_1+\be_2+\be_3|}}{4\al_1!\be_1!\cdots\be_3!}\\
\times (\zeta_{+}\hte\phi_1)^{(\al_1+\al_2)}_{(\be_1+\be_2)}({\tilde e}_1a^2\Gamma)^{(\be_1+\al_3)}_{(\al_1+\be_3)}(\zeta_{+}\hte\phi_1)^{(\be_2+\be_3)}_{(\al_2+\al_3)}+R
\end{align*}
where the sum is taken over $|\al_1+\be_1+\cdots+\be_3|=2$ and $R\in \mu^2 S(\la^{1+0},\bg)$ which follows from Lemma \ref{lem:toyama}. Here it can be checked that the second term is written as 
\[
c_1\zeta^2w^{-1}\la \hte^2+c_2\zeta w^{-1}\hte\phi_1+c_3w^{-1}\hphi_1^2\la+c_4w^{-1/2}\phi_1+c_5w^{-1/2}\zeta \hte \la
\]
with $c_i\in \mu^2 S(\la_{\mu}^{+0},\bg)$ modulo $\mu^2 S(w^{-1}\la_{\mu}^{+0},\bg)$. To estimate the first term let us write $c_1\zeta^2w^{-1}\la\hte^2={\mathsf{Re}}(c_1\zeta w^{-3/2}\hte\#\zeta w^{1/2}\la\hte)+R$ with $R\in \mu^2 S(\la^{1+0},\bg)$.  Then one can estimate $|{\mathsf{Re}}(c_1\zeta^2w^{-1}\hte^2 u,u)|$ by
\[
C\mu^3 n\|\zeta w^{-3/2}\hte u\|^2+C\mu n^{-1}\|\zeta w^{1/2}\la\hte
u\|^2+C\|\la_{\mu}^{+0}O(E)u\|^2.
\]
It is easy to see that $|((c_2\zeta w^{-1}\hte\phi_1+c_3w^{-1}\hphi_1^2\la+c_4w^{-1/2}\phi_1+c_5w^{-1/2}\zeta \hte \la
)u,u)|$ is bounded by $
C\ga^{-1/2}(\|\zeta w^{-3/2}\hte u\|^2+\|\zeta w^{1/2}\hte \la u\|^2)+C\ga^{1/2}\|O(E)u\|^2$. To end
the proof it suffices to apply the Fefferman-Phong inequality to obtain 
\begin{align*}
{\mathsf{Re}}({\tilde e}_1a^2\Gamma(\zeta_{+}\hte\phi_1u), \zeta_{+}\hte\phi_1u)\geq c\mu{\mathsf{Re}}(\Gamma(\zeta_{+}\hte\phi_1u), \zeta_{+}\hte\phi_1u)
-C\|\hte \phi_1u\|^2
\end{align*}
because  ${\tilde e}_1a^2-c\mu\geq 0$ with some $c>0$.
\qed
 
Similar arguments proving Lemma \ref{lem:foo} shows the estimate
\begin{align*}
{\mathsf{Re}}(a^2\chi_2\hphi_1^4\la^2 u,{\tilde e}_1\zeta_{+}^2\hte\Gamma u)\geq -C\ga^{-1/2}(\|\zeta_{+}w^{-3/2}\hte u\|^2+\|\zeta_{+}w^{1/2}\hte\la u\|^2)\\
-C\ga^{1/2}(\|\zeta_{+}w^{-1}\hte u\|^2+\|O(E)u\|^2).
\end{align*}
We turn to consider 
\begin{equation}
\label{eq:kaze}
((a_{1}\mu w^{-2}\hte+a_2w^{-1/2}\hte\phi_1
+a_3\phi_1)u,({\mathsf{Im}}{\tilde \la}) u).
\end{equation}
To handle \eqref{eq:kaze} we prepare a lemma.
\begin{lem}
\label{lem:nenga} We have 
\[
{\mathsf{Re}}(\Gamma u,v)\leq (\Gamma v,v)+(\Gamma w,w)+C(\|\la_{\mu}^{+0}v\|^2+\|\la_{\mu}^{+0}w\|^2).
\]
\end{lem}
%
\noindent
Proof:
Since $0\leq \Gamma\in S(\la_{\mu}^{1/2+0}, \bg)$ it follows from the Fefferman-Phong inequality that that $
(\Gamma u,u)\geq -C\|\la_{\mu}^{+0}u\|^2$ 
with some $C>0$. Thus with $L=\Gamma+C\la_{\mu}^{+0}$ we have $
(Lu,u)\geq 0$ so that $
|{\mathsf{Re}}(L u,v)|\leq (L u,u)+(L v,v)$ 
which proves the assertion. 
\qed
 
Write ${\mathsf{Re}}\,({\tilde e}_1\zeta_{+}^2\hte\Gamma\#a_1\mu w^{-2}\hte)=\mu{\mathsf{Re}}\,(\Gamma \zeta_{+}\hte w\la\#a\zeta_{+}w^{-1}\hte)+R$ with $R\in S(w^{-2},\bg)$ and apply  Lemma \ref{lem:nenga} to get
\begin{align*}
|{\mathsf{Re}}\,({\tilde e}_1\Gamma\zeta_{+}^2\hte u,a_1\mu w^{-2}\hte u)|\leq 
C\mu n^{-1}{\mathsf{Re}}\,(\Gamma \zeta_{+}\hte w\la u,\zeta_{+}\hte w\la u)\\
+C\mu n{\mathsf{Re}}\,(\Gamma a\zeta_{+}w^{-1}\hte u, a\zeta_{+}w^{-1}\hte u)\\
+C\mu(\|w\la u\|^2+\mu^{2}\|\zeta_{+}w^{-1}\hte u\|^2+\|w^{-1}\la_{\mu}^{+0}u\|^2).
\end{align*}
Since  $|{\mathsf{Re}}\,(\Gamma a\zeta_{+}w^{-1}\hte u, a\zeta_{+}w^{-1}\hte u)|\leq C\mu^{2}(\|\zeta_{+}w^{-3/2}\hte u\|^2+\|w^{-1}\la_{\mu}^{+0}u\|^2)$ we conclude
\begin{align*}
|{\mathsf{Re}}\,({\tilde e}_1\Gamma\zeta_{+}^2\hte u,a_1\mu w^{-2}\hte u)|\leq 
C\mu n^{-1}{\mathsf{Re}}\,(\Gamma \zeta_{+}\hte w\la u,\zeta_{+}\hte w\la u)
\leq C\|w\la u\|^2\\+C\mu^3 n\|\zeta_{+}w^{-3/2}\hte u\|^2+C(\|\zeta_{+}w^{-1}\hte u\|^2+\|w^{-1}\la_{\mu}^{+0}u\|^2).
\end{align*}
Similar arguments shows
\begin{align*}
|{\mathsf{Re}}\,({\tilde e}_1\Gamma\zeta_{+}^2\hte u,a_3\phi_1 u)|\leq 
\gamma^{-1/2}{\mathsf{Re}}\,(\Gamma \zeta_{+}\hte w\la u,\zeta_{+}\hte w\la u)\\
\leq C\|w\la u\|^2+C\gamma^{1/2}\mu^{-2}(\|w^{1/2}\phi_1 u\|^2+\|w^{-1}\la_{\mu}^{+0}u\|^2).
\end{align*}
Repeating similar arguments we conclude that \eqref{eq:kaze} is bounded by
\begin{align*}
C(\mu n^{-1}+\gamma^{-1/2}){\mathsf{Re}}\,(\Gamma \zeta_{+}\hte w\la u,\zeta_{+}\hte w\la u)+C\mu^2 n\|\zeta_{+}w^{-3/2}\hte u\|^2\\
+C\gamma^{1/2}(\|\zeta_{+}w^{-1}\hte  u\|^2+\|O(E)u\|^2).
\end{align*}
We finally consider the term $(q u, bu)$ with $b\in S(\la_{\mu}^{+0},\bg)$. Noticing ${\tilde\zeta}'\hte^{1/2}\in S(w^{1/2},\bg)$ one sees 
\begin{align*}
{\mathsf{Re}}(b\#a^2{\tilde \zeta}^2\hte\phi_1^2)={\mathsf{Re}}(ba{\tilde \zeta}\hte^{1/2}\phi_1\#a{\tilde \zeta}\hte^{1/2}\phi_1)+O(E)\cdot O(E)+O(E)
\end{align*}
and hence one obtains $|(qu,bu)|\leq C\|O(E)u\|^2$. We summarize
\begin{pro}
\label{pro:hagaki} There is $c>0$ and one can find $\ga_0>0$, $\mu_0>0$, $n_0>0$, $\nu_0>0$ such that we have
\begin{align*}
C\big\{\gamma((q+\mu\la)u,u)+\gamma^3\|u\|^2+\gamma\|{\tilde \Lambda}u\|^2
+{\mathsf{Re}}\,(({\mathsf{Re}}\, {\tilde Q}-T_1+{\bar\kappa}\mu\la)u,{\mathsf{Im}}\,{\tilde \lambda}\,u)\\
+\mu n\|\chi\zeta w^{-1/2}\tLa u\|^2\big\}
\geq cn \nu\mu (\Gamma(\zeta_{+} \hte\phi_1)u,\zeta_{+} \hte\phi_1 u)\\
+cn\mu|(\zeta^2\hte \rho \la^2 u,u)|+c n\mu\|\zeta\hte^{1/2}\rho^{1/2}\la  u\|^2
\end{align*}
for $\ga\geq \ga_0$, $0<\mu<\mu_0$, $n\geq n_0$ and $0<\nu\leq\nu_0$.
\end{pro}
%

\section{Estimates of error terms}

In this section we disregard error terms which are bounded by $\ga^2\|\la_{\mu}^{+0}u\|^2$ because we have $\ga^3\|\la_{\mu}^{3\varepsilon}u\|^2$ in \eqref{eq:minore}. 
We estimate ${\mathsf{Re}}({\tilde\La}u,({\mathsf{Im}}\,{\tilde Q}-T_2)u)$.
Recall
\begin{align*}
{\mathsf{Im}}\,{\tilde Q}-T_2=-\nu n\zeta_{+}^2a^2\hte\phi_1 \Gamma \{\phi_1,\hphi_2\}+a_1\mu w^{-3/2}\hte
+a_2w^{1/2}\hte\la \\
+a_3\hte\phi_1
+c_1O(E)+c_2O(E)+\mu S(\la_{\mu}^{+0},\bg)O(E).
\end{align*}
Thanks to Lemma \ref{lem:toyama} one can write 
\[
a^2\zeta_{+}^2\hte\phi_1\{\phi_1,\hphi_2\} \Gamma =\mu \zeta_{+}\# {\hat a} \Gamma \#\zeta_{+}\hte\phi_1+c_1\zeta_{+}\hte\la+c_2\zeta_{+}\phi_1
\]
with ${\hat a}=\mu^{-1}a^2\{\phi_1,\hphi_2\}$ modulo $\la_{\mu}^{+0}O(E)$ where $c_i\in \mu S(\la_{\mu}^{+0},\bg)$. Noticing $\zeta\zeta_{+}=\zeta_{+}$ from Lemma \ref{lem:nenga} it follows that
\begin{eqnarray*}
\nu n{\mathsf{Re}}(a^2\zeta_{+}^2\hte\phi_1\{\phi_1,\hphi_2\}
\Gamma u,\tLa u)
\leq \ep^{-1}n\nu^2\mu({\hat a}\Gamma(\zeta_{+}\hte \phi_1) u,(\zeta_{+}\hte \phi_1)u)\\
+\ep n\mu({\hat a}\Gamma \zeta(\tLa u),\zeta(\tLa u))+cn\nu\mu\|\zeta_{+}\rho^{-1/2}\tLa u\|^2
+cn\nu\mu\|\zeta_{+} \rho^{1/2}\hte \la u\|^2\\
+C(\|\la_{\mu}^{+0}\hte\phi_1 u\|^2+\|\la_{\mu}^{+0}\tLa u\|^2+\|\la_{\mu}^{+0}O(E)u\|^2)
\end{eqnarray*}
where $\ep>0$ will be determined later. We turn to estimate  
\[
((a_1\mu w^{-3/2}\hte
+a_2w^{1/2}\hte\la 
+a_3\hte\phi_1)u,{\tLa}u).
\]
It is easy to see that this is bounded by
\begin{align*}
C\ga^{-1/2}(\|\zeta w^{-3/2}\hte u\|^2+\|\zeta w^{1/2}\hte\la u\|^2)\\
+C\gamma^{1/2}(\|\tLa u\|^2
+\|\hte\phi_1 u\|^2
+\|\la_{\mu}^{+0}O(E)u\|^2).
\end{align*}
Finally we consider $|(c_1O(E)u+c_2O(E)u, \tLa u)|$. Recalling Corollary \ref{cor:QIm} it is easily seen that this term is estimated by 
\[
C\ga^{-1/2}(\|\zeta\chi w^{-1/2}\tLa u\|^2+\|\zeta\rho^{-1/2}\tLa u\|^2)+C\ga^{1/2}\|\la_{\mu}^{+0}O(E)u\|^2.
\]
Noting $\|w^{1/2}\phi_1u\|^2+\|\omega^{1/2}\phi_1u\|^2\geq \|\rho^{1/2}\phi_1 u\|^2-C\|O(E)u\|^2$ we summarize
\begin{pro}
\label{pro:amari:a} The term $|{\mathsf{Re}}(\tLa u,({\mathsf{Im}}\,{\tilde Q}-T_2)u)|$ is bounded by
\begin{align*}
 c\epsilon^{-1}n\nu^2\mu(\Gamma(\zeta_{+}\hte\phi_1)u,(\zeta_{+}\hte\phi_1)u)
+c\epsilon n\mu(\Gamma\zeta(\tLa u),\zeta(\tLa u))\\
+(cn\nu\mu+C\ga^{-1/2})(\|\zeta\rho^{-1/2}\tLa u\|^2
+\|\zeta \rho^{1/2}\hte \la u\|^2)\\
+C\ga^{-1/2}(\|\zeta w^{-3/2}\hte u\|^2
+\|\zeta\chi w^{-1/2}\tLa u\|^2)\\
+C\ga^{1/2}(\|\hte\phi_1 u\|^2+\|\tLa u\|^2+\|\la_{\mu}^{+0}O(E)u\|^2)
\end{align*}
where $c>0$ is independent of $\epsilon$, $\nu$, $\mu$ and $\ga$.
\end{pro}

We turn to consider the commutator $([D_0-{\mathsf{Re}}{\tilde \la},{\mathsf{Re}}\,{\tilde Q}-T_1]u,u)$. Recall
\[
\xi_0-{\mathsf{Re}}\,{\tilde \la}=\xi_0-\phi_1+\psi+n (b_2\hte+b_3\hphi_1^2) w^{-1/2}+R_1
\]
where $\psi={\tilde\zeta}\hte\phi_1+\chi_2\hphi_1^3\la$ and $R\in S(\la_{\mu}^{+0},\bg)$ and
\begin{align*}
{\mathsf{Re}}\,{\tilde Q}-T_1=q+\zeta_{+}(a_{1}\mu w^{-2}\hte+a_2w^{-1/2}\hte\phi_1
+a_3\phi_1)+\mu S(\la_{\mu}^{+0},\bg)O(E)
\end{align*}
where $q=\phi_2^2+2{\tilde \zeta}a^2\hte\phi_1^2+2\chi_2a^2\hphi_1^4\la^2$. Let us study $(\phi_2\{\xi_0-\phi_1+\psi,\phi_2\}u,u)$. Taking \eqref{eq:bunkai} into account it suffices to estimate 
\[
\nu(c_1\zeta_{+}^2\hte\phi_2\la u,u),\;\;(c_2\phi_2\hphi_1^2\la u,u),\;\; (c_3\hte\phi_2\phi_1 u,u)
\]
where $c_j\in \mu S(1,\bg)$. Write $
\zeta_{+}^2\hte\phi_2\la=(1-\chi^2)\zeta_{+}^2 \hte\phi_2\la+\chi^2\zeta_{+}^2 \hte\phi_2\la$ and  consider $M\zeta_{+}^2 \hte\phi_1^2-(1-\chi^2)\zeta_{+}^2 \hte\phi_2\la$ with a large positive constant $M$. Note 
\begin{eqnarray*}
M^2\zeta_{+}^2 \hte\phi_1^2-(1-\chi^2)\zeta_{+}^2 \hte\phi_2\la
=(M\phi_1)^2F
\end{eqnarray*}
where $0\leq F=\zeta_{+}^2 \hte\big(1-(1-\chi^2)\hphi_2\hphi_1^{-2}/M^2\big)\in S(1,g)$. Writing  $
(M\phi_1)^2F={\mathsf{Re}}(M\phi_1\#F\#M\phi_1)+R$ 
with $R\in  S(w^{-2},\bg)$ we obtain from the Fefferman-Phong inequality that
\begin{align*}
M^2(\zeta_{+}^2 \hte\phi_1^2 u,u)\geq ((1-\chi^2)\zeta_{+}^2 \hte\phi_2\la u,u)-C\|O(E)u\|^2.
\end{align*}
Consider now $
2w\chi^2\zeta_{+}^2 \hte\la^2-\chi^2\zeta^2 \hte\phi_2\la=(w^{1/2}\la)^2F$
with $0\leq F=\chi^2\zeta_{+}^2 \hte(2-\hphi_2w^{-1})\in S(1,\bg)$. Since $
(w^{1/2}\la)^2F={\mathsf{Re}}(w^{1/2}\la\#F\#w^{1/2}\la)+R$ with $R\in \mu S(\la,\bg)$ from the Fefferman-Phong inequality one has
\[
2(w\chi^2\zeta_{+}^2 \hte\la^2 u,u)\geq (\chi^2\zeta_{+}^2 \hte\phi_2\la u,u)-C\|w^{1/2}\la^{1/2}u\|^2-C\|O(E)u\|^2.
\]
Here we note that $
w^{1/2}\la^{1/2}\#w^{1/2}\la^{1/2}=w\la+R$ with $R\in S(1,\bg)$ and hence 
\begin{align*}
2(w\chi^2\zeta^2_{+}\hte\la^2 u,u)\geq (\chi^2\zeta_{+}^2 \hte\phi_2\la u,u)-C\|O(E)u\|^2.
\end{align*}
It is easy to see $|(c_3\hte \phi_1\phi_2u,u)|\leq C(\|\hte\phi_1 u\|^2+\|O(E)u\|^2)$ then we summarize
\begin{lem}
\label{lem:MANA} We have
\begin{eqnarray*}
|(\{\xi_0-\phi_1+\psi,\phi_2^2\} u,u)|\leq C\nu \mu(\chi^2\zeta_{+}^2 \hte w\la^2 u,u)\\
+C(\zeta_{+}^2\hte\phi_1^2 u,u)
+C(\|\hte\phi_1 u\|^2+\|O(E)u\|^2).
\end{eqnarray*}
\end{lem}
We next consider $\{\xi_0-\phi_1+\psi,{\tilde \zeta}a^2\hte\phi_1^2\}$ which is
\begin{align*}
\{\xi_0-\phi_1,{\tilde \zeta}a^2\hte\phi_1^2\}
+\{{\tilde\zeta}\hte\phi_1,a^2\phi_1\}{\tilde \zeta}\hte\phi_1
+\{\chi_2\hphi_1^3\la,{\tilde \zeta}a^2\hte\phi_1^2\}.
\end{align*}
It follows that  $\{\chi_2\hphi_1^3\la,{\tilde\zeta}a^2\hte\phi_1^2\}=c_1\hphi_1^4\la^2$ and $\{{\tilde \zeta}\hte \phi_1,a^2\phi_1\}{\tilde\zeta}\hte \phi_1=c_2\hte \phi_1^2$ from Lemma \ref{lem:toyama}. Since $
\{\xi_0-\phi_1,{\tilde\zeta}a^2\hte\phi_1^2\}=c_1\hte\phi_1^2+c_2\hphi_1^2\la\phi_2+c_3\hte\phi_1 \phi_2+c_4\hphi_1^4\la^2
$ by \eqref{eq:global:assump}, \eqref{eq:global:assump:b} and Lemma \ref{lem:DGB} we get
\begin{align*}
\{\xi_0-\phi_1+\psi,{\tilde \zeta}a^2\hte\phi_1^2\}=c_1\hte\phi_1^2+c_2\hphi_1^2\la\phi_2+c_4\hte\phi_1\phi_2+c_5\hphi_1^4\la^2.
\end{align*}
We then consider $\{\xi_0-\phi_1+\psi,\chi_2 a^2\hphi_1^4\la^2\}$ which is
\begin{align*}
\{\xi_0-\phi_1,\chi_2a^2\hphi_1^4\la^2\}
+\{{\tilde\zeta}\hte\phi_1,\chi_2a^2\hphi_1^4\la^2\}
+\{\chi_2\hphi_1^3\la,a^2\hphi_1\la\}\chi_2\hphi_1^3\la.
\end{align*}
A repetition of similar arguments shows
\begin{align*}
\{\xi_0-\phi_1+\psi,\chi_2 a^2\hphi_1^4\la^2\}=c_1\hte\phi_1^2+c_2\hphi_1^4\la^2+\hphi_1^2\la\phi_2.
\end{align*}
Therefore $|(\{\xi_0-\phi_1+\psi,a^2{\tilde \zeta}\hte\phi_1^2+a^2\chi_2\hphi_1^4\la^2\}u,u)|$ is bounded by 
\[
C|(\hte\phi_1^2 u,u)|+C\|O(E)u\|^2.
\]
Denoting $\zeta_{+}a_j$ by $a_j$ we turn to check 
$\{\xi_0-\phi_1+\psi,a_2w^{-1/2}\hte\phi_1\}$ where $a_2\in S(\la_{\mu}^{+0},\bg)$ of which support is contained in ${\rm supp}\,\zeta_{+}$. Remarking Lemmas \ref{lem:toyama} and \ref{lem:innsi} it is easy to see that 
\[
\{\xi_0-\phi_1+\psi,a_2w^{-1/2}\hte\phi_1\}=c_0\zeta^2 w^{1/2}\la\hte\phi_1+c_1\mu \zeta w^{-1/2}\hte \la+c_2\mu \la
\]
with $c_j\in \mu S(\la_{\mu}^{+0},\bg)$. Write $c_0\zeta^2 w^{1/2}\la\hte \phi_1={\mathsf{Re}}(c_0\zeta \hte^{1/2}\phi_1\#\zeta w^{1/2}\hte^{1/2}\la)+R_1$ and $c_1\zeta w^{-1/2}\hte\la={\mathsf{Re}}(\zeta w^{-3/2}\hte\#c_1w\la)+R_2$ with $R_i\in \mu S(\la^{1+0},\bg)$ we obtain the following estimate
\begin{align*}
|(c_0\zeta^2w^{1/2}\la\hte\phi_1 u,u)|\leq C\ga^{-1/2}(\|\zeta w^{1/2}\hte^{1/2}\la u\|^2+\|\zeta w^{-3/2}\hte u\|^2)\\
+C\ga^{1/2}\|c_0\zeta \hte^{1/2}\phi_1 u\|^2+C\|O(E)u\|^2.
\end{align*}
In order to estimate $\{\xi_0-\phi_1+\psi,a_1w^{-2}\hte\}$ we need to look at $a_1$ more carefully. Since $(w\phi)^{-1}\in S(\la_{\mu},g)$ the main part of $\{F,\log{\phi}\}$ is $w^{-1}\{F,\hphi_2\}$ by \eqref{eq:soji}. Therefore noticing \eqref{eq:Fai} it is not difficult to see from the proof of Lemma \ref{lem:tanaka} that $a_1$ has the form
\begin{equation}
\label{eq:marugame}
f(\zeta_{+})^{k_1}(\zeta_{+}')^{k_2}(\chi)^{k_3}(\chi')^{k_4}\hphi_1^{\ell_1}\hphi_2^{\ell_2}w^{s_1}\omega^{s_2}\rho^{s_3}(\log{\phi})^{\epsilon}
\end{equation}
where $f\in S(1,g_0)$ and $k_i$, $\ell_i\in\N$ and $s_i\in\R$, $\epsilon=0$ or $1$ which verify
\[
s_1+s_2+s_3+\ell_1/2+\ell_2\geq 0
\]
so that this is in $S(\la_{\mu}^{+0},\bg)$. Here we examine that $\xi_0-\phi_1+\psi$ commutes better against such terms of the form \eqref{eq:marugame} than  against general symbol in $S(\la_{\mu}^{+0},\bg)$. 
\begin{lem}
\label{lem:pisa} Denote $\La=\xi_0-\phi_1+\psi$ then $\{\La,\hphi_1\}$, $\{\La,\hphi_2\}$ and $\{\La,\hte\}$ is a linear combination of $\hphi_1$, $\hphi_2$ and $\hte$ with $\mu S(1,\bg)$ coefficients. We denote this by $\{\La,\hphi_1\}=\mu S(1,g_0)O(\Sigma)$ and so on.
\end{lem}
%
\noindent
Proof:
It follows easily from \eqref{eq:global:assump} and \eqref{eq:global:assump:b} that $\{\xi_0-\phi_1,\hphi_1\}$, $\{\xi_0-\phi_1,\hphi_2\}$ and $\{\xi_0-\phi_1,\hte\}$ are $O(\Sigma)$. Write $\psi=({\tilde \zeta}\hte+\chi_2\hphi_1^2)\phi_1$ and note Lemma \ref{lem:toyama} then the desired assertion for $\{\psi,\hphi_1\}$, $\{\psi,\hphi_2\}$ and $\{\psi,\hte\}$ follows immediately.
\qed
\begin{cor}
\label{cor:fujikawa} One can write $\{\La,w^{-1}\}=\mu S(w^{-1},\bg)+\mu S(w^{-2},\bg)O(\Sigma)$ and $\{\La,\omega^{-1}\}=\mu S(\omega^{-1},\bg)+\mu S(\omega^{-3/2},\bg)O(\Sigma)$ and that $\{\La,\rho^{-1}\}=\mu S(w^{-1},\bg)+\mu S(w^{-2},\bg)O(\Sigma)+S(\omega^{-3/2},\bg)O(\Sigma)$. We have also $\{\La,\zeta\}=c_1w^{-1}\hte+c_2w^{-1/2}$ with $c_i\in \mu S(1,\bg)$ and the same holds for $\{\La,\chi\}$.
\end{cor} 
Let us consider $\{\La,a_1\}$ where $a_1$ has the form \eqref{eq:marugame} with $k_1+k_2\geq1$. Since $(\chi)^{k_3}(\chi')^{k_4}\hphi_1\in S(w^{1/2},\bg)$ it follows from Lemma \ref{lem:pisa} and Corollary \ref{cor:fujikawa} that $\{\La,a_1\}$ can be written as $
c_0w^{-1}\hte+c_1w^{-1/2}+c_2\omega^{-1/2}$ 
with $c_i\in \mu S(\la_{\mu}^{+0},\bg)$. Since $\omega^{-1/2}w^{-1/2}\in \mu^{-1/2}S(\la^{1/2},\bg)$ then applying Lemma \ref{lem:pisa} and Corollary \ref{cor:fujikawa} again to $\{\La,w^{-2}\hte\}$ we conclude that
\[
\mu\{\La, a_1w^{-2}\hte\}=c_0\mu^3w^{-3}\hte^2+c_1\mu^{5/2}w^{-3/2}\hte \la^{1/2}+c_2\mu^2\phi_1+O(E)
\]
where $c_i\in S(\la_{\mu}^{+0},\bg)$. Writing $c_0w^{-3}\hte^2={\mathsf{Re}}(c_0w^{-3/2}\hte\#w^{-3/2}\hte)+R$ with $R\in S(w^{-2},\bg)$ and recalling that the support of $c_i$ are contained in the support of $\zeta_{+}$ we obtain the following estimate
\begin{equation}
\label{eq:koe}
\begin{split}
|\mu(\{\xi_0-\phi_1+\psi,a_1w^{-2}\hte\}u,u)|\leq C(\mu^3+\ga^{-1/2})\|\zeta w^{-3/2}\hte u\|^2\\
+C\ga^{1/2}\|O(E)u\|^2.
\end{split}
\end{equation}
If we write $b_1\omega^{-2}\hte=(b_1\omega^{-2}w^2)w^{-2}\hte$ and $b_3w^{-1}\omega^{-1}\hte=(b_3\omega^{-1}w)w^{-2}\hte$ then $b_1\omega^{-2}w^2$ and $b_3\omega^{-1}w$ have the same form \eqref{eq:marugame} and therefore we have the same estimate \eqref{eq:koe} for $|(\{\xi_0-\phi_1+\psi,b_1\omega^{-2}\hte+b_3w^{-1}\omega^{-1}\hte\}u,u)|$. Since the estimate $|(\{\xi_0-\phi_1+\psi,b_4\phi_1\}u,u)|\leq C(\|\hte\phi_1 u\|^2+\|O(E)u\|^2)$ is easy we summarize
\begin{pro}
\label{pro:amari:b} We have
\begin{align*}
|([D_0-{\mathsf{Re}}{\tilde \la},{\mathsf{Re}}\,{\tilde Q}]u,u)|\leq c\nu \mu(\chi^2\zeta_{+}^2 \hte w\la^2 u,u)\\
+(c\mu^3+C\ga^{-1/2})\|\zeta w^{-3/2}\hte u\|^2
+C(\zeta_{+}^2\hte\phi_1^2 u,u)\\
+C\ga^{-1/2}(\|\zeta \rho^{1/2}\hte^{1/2}\la u\|^2
+\|\zeta w^{-3/2}\hte u\|^2)
\\
+C\ga^{1/2}(\|\hte\phi_1 u\|^2+\|\zeta\hte^{1/2}\phi_1 u\|^2+\|O(E)u\|^2)
\end{align*}
where $c>0$ is indepensent of $\nu$, $\mu$ and $\ga$.
\end{pro}
%
\section{Lower order terms}

We finally handle the lower order terms. 
By \eqref{eq:ToneT} one can write
\[
T_2=\mu  c_0\hte\la+b_0\hte\phi_1+b_1\hphi_1^2\la+b_2\phi_2+b_3w^{1/2}\phi_1
\]
with $b_j\in \mu S(1,g)$ where $c_0=0$ for $\hte<0$ by assumption. Write $c_0\hte\la=c_0\zeta^2_{+}\hte\la+(1-\zeta^2_{+})c_0\hte\la$ where it is clear that we can write $(1-\zeta^2_{+})c_0\hte\la=b_4w\la$. We examine that one can write
\begin{align*}
(1-\chi^2)\zeta^2_{+}\hte\la=\om^{1/2}\rho^{-1}\zeta_{+}\#\rho\om^{-1/2}(1-\chi^2)\zeta_{+}\hte\la\\
+\om^{1/2}\rho^{-1}\zeta_{+}\#c\,\rho\om^{-1/2}\hte\la+R
\end{align*}
with $c\in S(\la_{\mu}^{-1/4},\bg)$ and $R\in S(\la^{1/2},\bg)$. Moreover ${\rm supp}\,c\subset {\rm supp}(1-\chi^2)$. Indeed since $\omega^{\pm1/2}\rho^{\mp1}\in S(\omega^{\pm1/2}\rho^{\mp 1},{\tilde g})$ then $\om^{1/2}\rho^{-1}\zeta_{+}\#\rho\om^{-1/2}(1-\chi^2)\zeta_{+}\hte\la$ can be written as $c\zeta_{+}\hte\la+R$ with $c\in S(\la_{\mu}^{-1/4},\bg)$ and $R\in \mu^{1/2}S(\la^{1/2},\bg)$. Write $c\zeta_{+}\hte\la=\om^{1/2}\rho^{-1}\zeta_{+}\#c\rho\om^{-1/2}\hte\la+R$ again we get the desired assertion. This proves
\begin{eqnarray*}
|((1-\chi^2)\zeta^2_{+}\hte\la u,\tLa u)|\leq C\ga^{-1/2}\|\om^{1/2}\rho^{-1}(\zeta_{+}\tLa u)\|^2
\\
+C\ga^{1/2}(\|c\rho\om^{-1/2}\hte\la u\|^2)
+\|\tLa u\|^2+\|O(E)u\|^2)
\end{eqnarray*}
with $c\in S(\la_{\mu}^{-1/4},\bg)$ where ${\rm supp}\, c\subset {\rm supp}(1-\chi^2)$. Now consider $
\|c\rho\om^{-1/2}\hte\la u\|^2$. Note that $c\rho \om^{-1/2}\in S(\om^{1/2},\bg)$ because if $c\neq 0$ then we have $C\hphi_1^2\geq w$ and hence $C^2\om^2\geq w^2\geq \hphi_2^2$. Thus it is clear $\om^2\leq \hphi_2^2+\om^2=\rho^2\leq (C^2+1)\om^2$ so that $\om^{1/2} \leq \rho\om^{-1/2}\leq (1+C')\om^{1/2}$. Hence it is easily seen that $c\rho\om^{-1/2}\hte\la\#c\rho\om^{-1/2}\hte\la=a\om \la^{3/2}+R$
with $a\in S(1,\bg)$ and $R\in S(\la,\bg)$ so that $
\|c\rho w^{-1/2} \hte\la u\|^2\leq C\|O(E)u\|^2$. 
We summarize
\begin{align*}
|((1-\chi^2)\zeta^2_{+}\hte\la u,\tLa u)|\leq C\ga^{-1/2}\|\rho^{-1}\om^{1/2}\zeta_{+}(\tLa u)\|^2\\
+C\ga^{1/2}(\|\tLa u\|^2+\|O(E)u\|^2).
\end{align*}
We turn to study $(\chi^2\zeta^2_{+}\hte\la u,\tLa u)$. Let us write $
\chi^2\zeta^2_{+}\hte\la=\chi \zeta_{+} w^{-1/2}\#\chi \zeta_{+}w^{1/2}\hte\la+c\hte\la+R$ 
with $c\in S(\la_{\mu}^{-1/2},\bg)$ and $R\in S(\la^{1/2},\bg)$ and hence we have
\begin{eqnarray*}
|(\chi^2\zeta^2_{+}\hte\la u,\tLa u)|\leq cn^{1/2}\|\chi \zeta_{+} w^{-1/2}\tLa u\|^2\\
+cn^{-1/2}\|\zeta_{+}\chi w^{1/2}\hte\la u\|^2
+C(\|O(E)u\|^2+\|\tLa u\|^2).
\end{eqnarray*}
Since it is clear that $|((b_0\hte\phi_1+b_1\hphi_1^2\la+b_2\phi_2+b_3w^{1/2}\phi_1)u,\tLa u)|$ is bounded by $C(\|\tLa u\|^2+\|\hte\phi_1 u\|^2+\|O(E)u\|^2)$ we get
\begin{pro}
\label{pro:LAMM}
We have
\begin{eqnarray*}
|(T_2 u,\tLa u)|\leq (c\mu n^{1/2}+C\ga^{-1/2})\|\chi \zeta_{+}w^{-1/2}\tLa u\|^2\\
+c\mu n^{-1/2}\|\chi\zeta_{+}w^{1/2}\hte\la u\|^2+C\ga^{-1/2}\|\rho^{-1}\om^{1/2}\zeta\tLa u\|^2\\
+C\ga^{1/2}(\|\tLa u\|^2+\|\la_{\mu}^{+0}O(E)u\|^2)
\end{eqnarray*}
with $c>0$ independent of $n$, $\nu$, $\mu$ and $\lambda$.
\end{pro}

We turn to consider $((T_1-{\bar\kappa}\mu\la)u, u)$.  
From Lemma \ref{lem:ameagari} it follows that  $\zeta_{-}^2h|\hte|\phi_1^2=h^{1/2}\phi_1\zeta_{-}|\hte|^{1/2}\#h^{1/2}\phi_1\zeta_{-}|\hte|^{1/2}+R$ with $R\in \mu^2 S(w^{-2},g)$. 
Noting $h^{1/2}\zeta_{-}|\hte|^{1/2}\in S(w^{1/2},g)$ we see 
\[
\phi_2\#h^{1/2}\phi_1\zeta_{-}|\hte|^{1/2}-h^{1/2}\phi_1\zeta_{-}|\hte|^{1/2}\#\phi_2=\{\phi_2,h^{1/2}\phi_1\zeta_{-}|\hte|^{1/2}\}/i+R
\]
with $R\in \mu^2 S(w^{-3/2},g)$. Here since $h=\mu{\hat c}\{\phi_1,\hphi_2\}^{-1}$  we have
\[
\{h^{1/2}\phi_1\zeta_{-}|\hte|^{1/2},\phi_2\}=\mu^{1/2}\zeta_{-}({\hat c}\{\phi_1,\hphi_2\}|\hte|)^{1/2}e+c\mu w^{1/2}\phi_1
\]
with $c\in S(1,g)$ thanks to Lemma \ref{lem:ameagari} because $\phi_2^{(\al)}\in \mu S(w,g)$ for $|\al|=1$. Then the following estimate  follows easily
\[
\mu^{1/2}\big(\zeta_{-}({\hat c}\{\phi_1,\hphi_2\}|\hte|)^{1/2}e\, u,u\big)\leq (\phi_2^2 u,u)+(h\zeta_{-}^2|\hte| \phi_1^2 u,u)+C\mu^{1/4}\|O(E)u\|^2
\]
because $\|w^{-1/2}u\|^2\leq C\mu^{-1/2}\|\la^{1/2}u\|^2$. From \eqref{eq:mao} it follows that
\begin{equation}
\label{eq:kouti}
\mu^{1/2}\zeta_{-}({\hat c}\{\phi_1,\hphi_2\}|\hte|)^{1/2}\,e+T_1\geq 2{\bar\kappa} \mu \la-C\mu w^{1/2}\la
\end{equation}
with some $C>0$. In fact if $\hte\leq -b_3w$ then $\zeta_{-}=1$ and the assertion by \eqref{eq:mao}. If $-b_3w\leq \hte\leq 0$ then we have $Cw^{1/2}\la\geq \mu^{-1/2}({\hat c}\{\phi_1,\hphi_2\}|\hte|)^{1/2}\,e$ and hence the assertion.  Since $S(\la,G)\subset S^1_{1,1/2}$ the Fefferman-Phong inequality gives 
\[
\mu^{1/2}({\hat c}\zeta_{-}(\{\phi_1,\hphi_2\}|\te|)^{1/2}e\, u,u)+(T_1 u,u)\geq 2{\bar\kappa} \mu(\la u,u)
-C\mu\|O(E)u\|^2.
\]
We summarize
\begin{pro}
\label{pro:yazu} We have
\begin{align*}
(\phi_2^2 u,u)+(h\zeta_{-}^2|\hte| \phi_1^2 u,u)+((T_1-{\bar\kappa}\mu\la)u,u)\\
\geq {\bar\kappa} \mu(\la u,u)
-C\mu^{1/4}\|O(E)u\|^2.
\end{align*}
Similarly $
(\la_{\mu}^{2\varepsilon}\phi_2^2 u,u)+(\la_{\mu}^{2\varepsilon}h\zeta_{-}^2|\hte| \phi_1^2 u,u)+(\la_{\mu}^{2\varepsilon}(T_1-{\bar\kappa}\mu\la)u,u)$ is bounded from below by ${\bar\kappa} \mu(\la_{\mu}^{2\varepsilon}\la u,u)
-C\mu^{1/4}\|\la_{\mu}^{\varepsilon}O(E)u\|^2$.
\end{pro}
Finally we estimate ${\mathsf{Re}}((T_1-{\bar\kappa}\mu\la)u,({\mathsf{Im}}{\tilde \lambda})\,u)$.  Since $\zeta_{-}\zeta_{+}=0$ then from \eqref{eq:kouti} we see that ${\mathsf{Re}}((T_1-{\bar\kappa}\mu\la)u,({\mathsf{Im}}{\tilde \lambda})\,u)$ is bounded from below by
\[
 {\mathsf{Re}}(({\bar\kappa}\mu\la-C\mu w^{1/2}\la)u,{\tilde e}_1\Gamma\zeta_{+}\hte  u)-C\|\la_{\mu}^{+0}O(E)u\|^2.
\]
Note that ${\mathsf{Re}}({\tilde e}_1\Gamma\zeta_{+}\hte \#({\bar\kappa}\mu\la-C\mu w^{1/2}\la)={\bar\kappa}\mu {\tilde e}_1\Gamma\zeta_{+}\hte\la+c\zeta_{+}w^{-1/2}\hte\la+R$ with $c\in S(\la_{\mu}^{+0},\bg)$ and $R\in \mu S(w^{-1},\bg)$. Since $0\leq {\tilde e}_1\Gamma\zeta_{+}\hte \la\in S(w^{-1}\la^{1+0},\bg)$ and noting $c\zeta_{+}w^{-1/2}\hte\la={\mathsf{Re}}(\zeta_{+}w^{-3/2}\hte\#cw\la)+R$ with $R\in S(\la,\bg)$ one can see that ${\mathsf{Re}}({\tilde e}_1\Gamma\zeta_{+}\hte \#({\bar\kappa}\mu\la-C\mu w^{1/2}\la)$ has a bound from below $-
C(\ga^{-1/2}\|\zeta_{+}w^{-3/2}\hte u\|^2+\ga^{1/2}\|\la_{\mu}^{+0}O(E)u\|^2)$. 
Therefore we obtain
\begin{lem}
\label{lem:simekiri} We have
\[
{\mathsf{Re}}((T_1-{\bar\kappa}\mu\la)u,({\mathsf{Im}}{\tilde \lambda})\,u)\geq -C\ga^{-1/2}\|\zeta_{+}w^{-3/2}\hte u\|^2-C\ga^{1/2}\|\la_{\mu}^{+0}O(E)u\|^2.
\]
\end{lem}
We first choose $\epsilon>0$ small so that 
$c\epsilon n\mu (\Gamma(\zeta(\tLa u),\zeta(\tLa u))$ in Proposition \ref{pro:amari:a} can be controlled by the corresponding  term in Proposition \ref{pro:akasaka}. We next choose $\nu>0$ small so that $cn\nu\mu\|\zeta\rho^{-1/2}\tLa u\|^2$ and
\[
 c\epsilon^{-1}n\nu^2\mu(\Gamma(\zeta_{+}\hte\phi_1)u,(\zeta_{+}\hte\phi_1)u)
+cn\nu\mu\|\zeta \rho^{1/2}\hte \la u\|^2
\]
in Proposition \ref{pro:amari:a} will be small against the corresponding  terms in Propositions \ref{pro:akasaka} and \ref{pro:hagaki}. We then choose $n$ such that $\mu^3\|\zeta w^{-3/2}\hte u\|^2$ in Proposition \ref{pro:amari:b} can be controlled by Proposition \ref{pro:wwLa} and $c\mu n^{1/2}\|\chi\zeta_{+}w^{-1/2}\tLa u\|^2+c\mu n^{-1/2}\|\chi\zeta_{+}w^{1/2}\hte\la u\|^2$ in Proposition \ref{pro:LAMM} can be estimated by Propositions \ref{pro:akasaka} and \ref{pro:hagaki}.  Finally we choose $\mu>0$ small enough and then $\ga>0$ enough large so that $\mu n^4$ to be small and $\ga\mu^4$ to be large. Then combining Propositions \ref{pro:akasaka}, \ref{pro:wwLa}, \ref{pro:wwLam}, \ref{pro:kasitu}, \ref{pro:hagaki}, \ref{pro:amari:a}, \ref{pro:amari:b}, \ref{pro:LAMM} and \ref{pro:yazu} we obtain a desired weighted energy estimates. Once we obtain energy estimates in order to conclude the well-posedness of the Cauchy problem it suffices to apply \cite[Theorem 1.1]{Ni5}.


\end{document}